\documentclass[12pt]{article}
\usepackage{amsmath,amsfonts, amssymb,verbatim}
\newtheorem{pkt}{}[section]
\newcommand{\bpk}{\begin{pkt}\rm }
\newcommand{\epk}{\end{pkt}}

\newcommand{\QQ}{{^*\mathbb{Q}}}

\newcommand{\CC}{{^*\mathbb{C}}}

\newcommand{\RR}{{'\R'}}
\newcommand{\inv}{^{-1}}
 
\newcommand{\G}{\mathrm{G}}
\newcommand{\U}{\mathbb{U}}
\newcommand{\uu}{{\bf u}}
\newcommand{\vv}{{\bf v}}
\newcommand{\ww}{{\bf w}}
\newcommand{\uuu}{\mathfrak{u}}
\newcommand{\vvv}{\mathfrak{v}}

\newcommand{\Oo}{\mathrm{O}}

\newcommand{\ii}{\mathbf{i}}
\newcommand{\jj}{\mathbf{j}}
\newcommand{\Dom}{\mathbb{V}}
\newcommand{\Pp}{\mathrm{P}}
\newcommand{\Qq}{\mathrm{Q}}
\newcommand{\GG}{{\cal G}}

\newcommand{\R}{{\mathbb R}}
\newcommand{\Q}{{\mathbb Q}}
\newcommand{\Z}{{\mathbb Z}}
\newcommand{\N}{{\mathbb N}}
\newcommand{\C}{{\mathbb C}}
\newcommand{\A}{\mathcal{A}}
\newcommand{\B}{\mathcal{B}}
\newcommand{\be}{\begin{equation}}
\newcommand{\ee}{\end{equation}}

\newcommand{\LL}{{\cal L}}
\newcommand{\HH}{\mathcal{H}}
\newcommand{\Ss}{\mathbb{S}}
\newcommand{\mfr}{\mathfrak{m}}

\newcommand{\Sss}{'\mathbb{S}' }

\newcommand{\K}{{\rm K}}

\newcommand{\acl}{{\rm acl}}
\newcommand{\NN}{\mathcal{N}}

\newcommand{\ra}{\rangle}
\newcommand{\la}{\langle}

\newcommand{\F}{{\rm F}}

\newcommand{\xx}{{\bf x}}

\newcommand{\s}{\mathbf{s}}
\newcommand{\pp}{\mathfrak{p}}
\newcommand{\lfr}{\mathfrak{l}}
\newcommand{\ZZ}{{^*\mathbb{Z}}}

\newcommand{\lm}{\mathsf{lm}}

\newcommand{\e}{\mathrm{e}}

\newcommand{\cc}{\mathbf{c}}

\title{On the logical structure of physics}
\author{Boris Zilber}

\begin{document}
\maketitle
\abstract{ One of the main claims of the paper is that  Dirac's calculus and broader theories of physics can be treated as theories written in the language of Continuous Logic. Establishing its true interpretation (model) is a model theory problem.  
The paper introduces such a model for the fragment which covers ``free theories'', that is physical theories with Gaussian (quadratic) potential.
  The model is  pseudo-finite (equivalently, a limit of finite models), based on a pseudo-finite field in place of the field of complex numbers. 
  The advantage of this unusual  setting is that it treats the quantum and the statistical mechanics as just domains  in the same model and explains Wick rotation as a natural transformation of the model corresponding to a shift in scales of physical units.  }
\tableofcontents

\section{Introduction}
\bpk 
The success of model theory in its numerous applications in various areas of mathematics is due, in the first place, to its focus on the fundamental questions: what is the adequate mathematical language of the area? what is the structure that is being studied? 

We start here by asking the same questions about physics or rather some parts of physics. 
More precisely our focus of interest is in foundations of quantum theory and of statistical physics, the two areas differing by the nature of physical processes and the scales of magnitudes,  using different mathematical platforms connected by the formal mathematical trick called ``the Wick rotation''.    

Statistical physics was properly established in the 19th centure when it was notices that the probability of a state (of a gas) with certain energy $E$ and temperatire $T$ is proportional to 
$\e^{-\frac{E}{T}}.$

Several decades later quantum physics started with the basic observation that the state of a quantum system of energy $E$ at time $t$ can only be adequately observed by attributing to it the complex number of norm 1,
$\e^{iEt},$
called the probability amplitude.


This settings generate certain logic which in the first case one can recognise as a form of the {\em probabilistic logic}. The second case is in fact more general and can be classified (as yet informally) as the
{\em continuous logic with values in $\C.$}  

The typical $n$-ary predicate in the theory has the form
\be \label{universe}\U^n \to \C; \ \ \bar{x}\cdot \uuu \mapsto c\cdot\e^{f(\bar{x})}\ee
where $\U$ is the domain of physical units, $\uuu\in \U,$ $f: \C^n\to \C$ some  function (typical for Hamiltonian mechanics) possibly with complex coefficients, $c$ a constant. 

The quantifiers are of the form
$$ \e^{f(\bar{x},y)} \mapsto \int_\R \e^{f(\bar{x},y)}\cdot \e^{a(y)}\ dy$$
and in basic cases the right-hand can be calculated to be   of the form $b\cdot \e^{g(\bar{x})}.$ 

In statistical theory $f$ and $a$ are real functions. In quantum theory these are typically  similar functions obtained by changing $f(x,y)$ and $a(y)$ to  $if(x,y)$ and $ia(y).$ The quantum mechanical version is in fact the calculus in Dirac's formalism.

The Wick rotation effect is that the calculation of quantifiers for real $f$ and $a$ return the same result $if(x,y)$ and $ia(y)$ with respective changes from     $b\cdot \e^{g(\bar{x})}$
to  $b'\cdot \e^{ig(\bar{x})},$ appropriately defined $b'.$  The effect can be explained mathematically in many cases but the physical nature of the formal link between statistical and quantum theory remains a mystery.
\epk
\bpk In this paper we take seriously the assumption that the logical setting of physics is that of continuous logic (CL) and treat the  laws of physics written in terms of  CL-formulas as {\em axioms of a CL-theory}. The problem is to find an {\em interpretation} of the axioms, that is a class of yet to be defined structures: continuous, finite, or pseudo-finite . 

(In \cite{Udi} E.Hrushovski solves a similar, albeit in fact inverse, problem of writing down a CL-theory for a class of finite and pseudo-finite structures.)

First of all we need to determine the {\bf universe} $\U$ as in (\ref{universe}) and the structure on it.

A crucial question of significance in physics is that whether $\U$ is discrete or continuous. In statistical mechanics $\U$ is discrete, and even  finite of huge size, by definition. In quantum theory it is currently universally accepted that $\U$ is discrete (which follows from the fact that distinct points of the space has to be  at the distance at least of the Planck length). It is also consistent to assume more generally the finiteness of the universe if we accept that the age of the universe is finite. These issues are being actively discussed in physics literature.

Now, if we accept that $\U$ is finite or pseudo-finite (in the model theory sense) then the range of the maps in (\ref{universe}) does not need to be a continuous field, in fact it makes sense to consider   logical values in a finite or pseudo-finite field $\F_\pp$ in place of $\C.$ However, there is an obvious difficulty of practical nature. 
R.Penrose in his book {\em The Road to Reality} writes about the prospect  of using finite fields in physics:

{\em ... It is unclear whether such things really have a significalnt role to play in physics, although the idea has been revived from time to time. If $\F_p$ were to take the place of the real-number system, in any significant sense, then $p$ has to be very large indeed. ... To my mind, a physical theory  which depends fundamentally upon some absurdly enormous prime number would be a far more complicated (and improbable) theory than one that is able to depend on a simple notion of infinity. Nevertheless, it is of some interest to pursue these matters. ...}

One of the main result of our paper is that the setting with pseudo-finite $\F_\pp$ is quite easily convertible to the setting over the field $\C.$ Conversely the continuous calculus over $\C$ (and $\R$) can be translated into some meaningful calculations over $\F_\pp$ without loss.

\epk
\bpk The mathematical tool which is behind the passage between the pseudo-finite setting and $\C$ is the following diagram first established in \cite{FP} and worked out in the current paper in more detail:
\be\label{Dgm}\begin{array}{lllllll}
\ \ \ \ \ \ \ \ \ \ \ \ \ \  \lm_\U\ \ \ \ \ \ \ \ \\
\ \ \ \ \ \U\ \ \ \ \ \longrightarrow\ \ \ {\bar{\C}}\\
\ \ \ \ \ \ \ |\ \ \ \ \ \ \ \ \ \ \ \ \ \ \ |\ \ \\          
\exp_\pp\ |\ \ \ \ \ \ \ \ \ \ \ \ \ \ \ |\exp\\ 
\ \ \ \ \ \ \ \downarrow\ \ \ \ \ \ \ \ \ \ \ \ \ \ \downarrow\ \ \\
\ \ \ \ \ \ \ \ \ \ \ \ \ \  \lm_\F\ \ \ \ \ \ \ \ \\
\ \ \ \ \ \ \F_\pp \ \ \ \  \longrightarrow \ \  \ \bar{\C}
\end{array}\ee
Here $\U$ is a pseudo-finite additive group and 
 $\exp_\pp$ is a surjective homomorphism of the additive group onto $\F_\pp^\times.$
The horizontal arrows are ``limit'' maps $\mathsf{lm}_\U$ and $\mathsf{lm}_\F$ respectively, $\mathsf{lm}_\F$ is a place from the field $\F_\pp$ onto $\bar{\C}=\C\cup \{ \infty\},$ and $\mathsf{lm}_\U$ is appropriate map for the additive structures.

It is crucial that the limit maps are rather well-controlled, in particular  certain natural multpilicative subgroups $'\Ss',\, \RR_+ \subset \F_\pp$ map onto the unit circle $\Ss\subset \bar{\C}$ and the non-negative reals $\R_+\subset \C$  which allows  to mimic polar coordinates of $\C$ in $\F_\pp$ and thus develop a working analogue of continuous complex 
calculus.
\epk
\bpk In essense, { the main  model theory result of the paper is the construction of a continuous-model ultraproduct of quite intricate finite structures interpreted in 
$(\U_{p}, \F_p)$ which include wave-functions,  finite-dimensional ``Hermitian and Euclidean Hilbert spaces'' over $\F_p,$ linear operators on the spaces and other relevant constructs.    } 

One of the main gains of the pseudo-finite setting in terms of foundations of physics is in explaining the effect of Wick rotation as the transformation/homomorphism of $\U$ caused by the multiplication
$$u\mapsto \ii\cdot u$$
where $\ii$ is a non-standard integer such that 
$$\lm_\U(\ii\cdot u)=i \cdot \lm_U(u)\ \ \ (\ i=\sqrt{-1}\ ).$$ 

The action of the ``huge'' integer $\ii$ shifts a subdomain $\Dom_\uuu$ of the universe $\U$ to the subdomain $\Dom_\vvv= \ii\cdot \Dom_\uuu.$
We associate $\Dom_\uuu$ with physics at the scale of statistical mechanics and the Euclidean Hilbert space formalism and $\Dom_\vvv$ with quantum mechanics and Hermitian Hilbert spaces. Thus the shift by $\ii$ is the mathematical form of the change of scales in physics which manifests itself as the Wick rotation. 
\epk
\section{Main results}

\bpk
We define $\U=\U_{\pp,\lfr}$ as the quotient of the additive group of non-standard integers $\ZZ$ by the ideal generated by the number $(\pp-1)\cdot\lfr,$
$$\U:= \ZZ/(\pp-1)\cdot\lfr$$
for  distinguished parameters: prime $\pp$ and a highly divisible $\lfr.$  $$\exp_\pp: \U\to \F_\pp^\times$$ is a homomorphism with the kernel $(\pp-1)\cdot \U$ and $\ii\in \ZZ$ divides $(\pp-1)$ and is divisible by $\lfr.$
 \epk
\bpk\label{LE} {\bf Logical evaluations of $\U$ in $\F_\pp.$}  We think of $\F_\pp$ 
 as 
 the {\bf domain of logical values} as opposed to  $\U$ and $\U^M$ as  {\bf domains of physical units}. We think of definable operations on $\F_\pp$ as logical connectives, where {\bf definable} means  interpretable in the non-standard model of arithemetic $\ZZ$  with distinguished parameters $\pp,$ $\lfr$ and the map $\exp_\pp.$ 
 

The map $\exp_\pp: \U\to \F_\pp$ gives rise to evaluations of ``physical models'' and leads to a basic notion of the Hilbert space formalism of physics, a state (over $\F_\pp$): 
 $$\varphi: \Dom\to \F_\pp$$ where 
$\Dom\subset \U^M$ is a specific subdomain. A basic state (basic predicate) $\varphi$ has the form $$\varphi(\bar{x})=\exp_\pp(f(\bar{x})\cdot \uuu)$$
where $f(\bar{x})$ is a polynomial over $\Z,$ $\bar{x}\in \ZZ^N$ and $\uuu\in \Dom.$ A general state is obtained from basic ones by using  ``logical connectives''.


The states (which are of course a certain kind of coordinate functions) form  an $\F_\pp$-linear  space $\HH_\Dom$ of pseudo-finite dimension, with natural choices of orthonormal bases and well-defined  inner product with values in $\F$ or a well-controlled extension of $\F.$ Moreover, we consider definable linear maps on $\HH_\Dom,$ analogues of linear unitary operators playing an important role in physics (such as the Fourier transform and time evolution operators). 

The definable family  of {\bf position states} 
$$\LL_\uu=\{ \uu_r:\ r\in \Dom\}; \ \ \uu_r(s):=\delta(r-s)$$  
(Kronecker-delta) forms a basis of $\HH_\Dom.$ A definable injective linear operators $\A$ give rise to other definable bases
$$\LL_{\A\uu}=\{ \A\uu_r:\ r\in \Dom\}.$$
 
  \medskip

 In regards to model-theoretic formalism, unlike the traditional approach (see e.g. \cite{HartCMT}), 
$\HH_\Dom$ is not considered to be a universe of a structure. Instead, we consider the multisorted structure on sorts $\LL_\psi,$ definable bases,  with linear maps between these, together with the $\F_\pp$-linear space $\HH_\Dom$ interpretable in the sorts. 
 
\epk
\bpk 
In the current paper we restrict the study to  so called {\bf Gaussian setting}. This means that the $f(\bar{x})$ defining basic states $\varphi$ are 
quadratic forms on $\ZZ$ (in fact on the ring $\K_\NN=\ZZ/\NN,$ for $\NN=|\Dom|$)
 and the operators are of the form
\be\label{Aintro}\A: \ \varphi\mapsto \frac{1}{\sqrt{\NN}}\sum_{r\in \K_\NN} \exp_\pp(a(q,r)\uuu)\cdot \varphi(r),\ee  
where $\uuu\in \Dom,$ an element of order $\NN$ and $a(q,r)$ is a quadratic form.
These are discrete analogues of quantum mechanics unitary operators for the free particle and the harmonic oscillator. We believe  free fields theories should be representable in the setting once we switch to considering domains $\Dom\subset \U^M$ for appropriate infinite pseudo-finite $M.$ 

\medskip

 We are in particular interested in studying relationship between two domains  $\Dom_\uuu$ and $\Dom_\vvv$ in $\U$ determined by the choices of units, $\uuu$ and $\vvv$ respectively, of different scales: $\vvv=\ii\cdot \uuu$,  for $\ii$ described above. 
 
 \epk
 \bpk Our first main result,  Thm \ref{Th1}, establishes for the discrete pseudo-finite model  a  unifying treatment of ``physics'' over
 the two domains: 
 
 \medskip
 
  In fact, $$\Dom_\vvv\subset \Dom_\uuu$$
 both equipped with additive non-standard metrics defined in terms of units $\vvv$ or $\uuu$ respectively. 
 
 The mulitplication by $\ii$ determines the projection
  $$\ii: \ \Dom_\uuu\twoheadrightarrow \Dom_\vvv=\ii \Dom_\uuu.$$
  
Given a  state $\varphi$ on $\Dom_\uuu,$ its restriction to $\Dom_\vvv$  is a state on $\Dom_\vvv,$ which we denote $\varphi^\ii.$  In fact, if
$\varphi=\exp_\pp( f(r)\uuu)$ then $\varphi^\ii=\exp_\pp( \ii f(r)\uuu).$

Respectively the action of a linear operator 
$\A$ on   $\varphi\in \HH_{\Dom_\uuu}$ becomes the action of some well-defined linear operator $\A^\ii$ on  $\varphi^\ii\in \HH_{\Dom_\vvv},$
$$\A^\ii\varphi^\ii=(\A \varphi)^\ii.$$ 

A formal inner product on the spaces transforms correspondingly
$$\la \varphi^\ii| \psi^\ii\ra=\la \varphi| \psi\ra^\ii,$$
where we consider both a formal-Euclidean and a formal-Hermitian versions of inner product. 

This gives us the {\bf isomorphism of the  structures}
\be\label{morphism1}\{ \}^\ii: \ \HH_{\Dom_\uuu}\to  \HH_{\Dom_\uuu}\ee

\medskip

Note that a tensor product powers  $\HH_{\Dom}^{\otimes M}$ for pseudo-finite $M$ is  interpretable
in $\HH_{\Dom}$ (or rather in the underlying structure $(\U;\F_\pp)$).
Thus the  picture can be generalised to pseudo-finite-dimensional setting 
$$\Dom_\vvv^M\subset \Dom_\uuu^M\subset \U^M$$  
with the Hilbert spaces replaced by tensor product powers  $\HH_{\Dom_\vvv}^{\otimes M}$ and $\HH_{\Dom_\uuu}^{\otimes M}$ respectively.
 \epk
 \bpk Our final task is to recast the pseudo-finite ultraproduct of finite structures underlying the formal Hilbert spaces  
 $\HH_{\Dom_\uuu},$ $\HH_{\Dom_\vvv}$ 
 as a continuous logic (CL) structure with values in the complex numbers. The key to the construction are the limit maps $\lm_\U$ and $\lm_\F$ shown in diagram (\ref{Dgm}).
  

First we have to convert the domains $\Dom_\uuu$ and $\Dom_\vvv$ into metric spaces 
presentable as countable unions of finite diameter subspaces. That is not obvious since the nonstandard-valued dimeters of    $\Dom_\uuu$ and $\Dom_\vvv$ (in units $\uuu$ and $\vvv$ respectively) are large nonstandard numbers, however we follow these metrics and define, for $\Dom:=\Dom_\uuu$ and  $\Dom:=\Dom_\vvv$ respectively,
$$ \Dom_{|n}=\{ u\in \Dom:\ \mathrm{dist}(0,u)\le n\},$$ $$ \Dom_{|n}^\lm:=\lm_\U (\Dom_{|n})\mbox{ and }\Dom^\lm:= \bigcup_{n\in \N} \Dom_{|n}^\lm $$

It turns out that
 $$\Dom^\lm_\uuu=\R\mbox{ and }\Dom^\lm_\vvv=\imath\R$$
 
For a state $\varphi$ as above we define
$$\varphi^\lm: \ \Dom^\lm \to \C;\ \ \varphi^\lm(r^\lm):= \lm_\F(\sqrt{\NN} \varphi(r))$$ 
Note that we have to use a normalising coefficient, an infinite pseudo-finite number $\sqrt{\NN}$ ($\NN$ is $\NN_\uuu$ or $\NN_\vvv,$ respectively) in order to produce a meaningful wave-functions $\R\to \C.$ As a result, inner products also have to be renormalised, which agrees with the Dirac delta-function renormalisation of respecive integral formulas. 

Crucially, in case of $\Dom_\uuu$ we find it technically necessary to choose the Euclidean inner product, and in case of $\Dom_\vvv$ the
Hermitian inner product.  

Finally, the linear operators $\A$ of the form (\ref{Aintro}) becomes the integral operator on the Hilbert space of functions $\phi,$
\be \label{Alm} \A^\lm: \phi\mapsto \int_{\R} \e^{\alpha(s,x)}\cdot \phi(x)\ dx.\ee
In terms of continuous logic, this is an existential quantifier (bounded by condition $\e^{\alpha(s,x)}$). 

\medskip

We refer to the two structures as 
$$\HH_\R:=\HH^\lm_{\Dom_\uuu}\mbox{ and }\HH_{\imath \R}:=\HH^\lm_{\Dom_\vvv}$$
with Euclidean and Hermitian inner products respectively.
\epk
\bpk \label{Final} The final result presented in section \ref{CL}  can be summarised as  the following  Theorem compairing the two continuous logic structures:   

\medskip



{\em The map (\ref{morphism1}) passes to a morphism of CL-structures
$$\{ \}^\ii:\ \HH_\R \to \HH_{\imath \R}$$
realised by the bijection on the CL-position states sorts
$$ \LL^\lm_\uu\to  \LL^\lm_{\uu^\ii}; \ \ \uu_r\mapsto \uu_{\imath r},\ \ r\in \R$$
and general ket-sorts
$$ \LL^\lm_\psi\to  \LL^\lm_{\psi^\ii}; \ \ \psi_r\mapsto \psi_{\imath r},\ \ r\in \R$$
which commute with the tranformation of the integral operators (\ref{Alm})
$$\int_\R \e^{\alpha(s,x)}\cdot \phi(x) \ dx \mapsto \int_\R \e^{\imath \alpha(s,x)}\cdot \phi^\ii(x) \ dx$$
and induces the $\R$-bilinear bijective maps between inner products
$$\la  \psi_r|\varphi_s\ra_\mathrm{E} \mapsto  \la  \psi^\ii_r|\varphi^\ii_s\ra_\mathrm{H};\ \ r,s\in \R$$  
(Euclidean on LHS and Hermitian on RHS)

This gives a full account on the Wick rotation in  Gaussian setting.

\medskip

Also, treating the integral operators as quantifiers, the structure allows quantifier elimination.}

\epk
\bpk {\bf Future directions.} The extension of the Gaussian setting to free fields theories seems quite feasible. The more general setting requiring perturbation methods is more challenging but does not seem impossible since tools of real and complex analysis are reasonably available in  $\F_\pp$-setting. 
\epk

\section{Definability and scales in $\U.$ }
\bpk\label{*C} As in \cite{FP},
let $$\CC=\C^P/D,\ \ \ZZ=\Z^P/D$$ be  ultrapowers of the field of complex numbers and of the ring of integers by a non-principal ultrafilter on the set of prime numbers $P.$

 Recall that by definition we have a representation of $\U_{\pp,\lfr}$ and $\F_\pp$ together with $\exp_\pp$ in $(\ZZ; +,\cdot, \pp,\ii,\lfr)$ as quotients
 $$\U_{\pp,\lfr}= \ZZ/(\pp-1)\lfr;\ \ \ \F_\pp=\ZZ/\pp$$
and
 $$\exp_\pp: \eta\cdot \hat{1}\mapsto \epsilon^\eta, \mbox{ for some } \epsilon\in \F_\pp\mbox{ and each }\eta\in \ZZ$$
 such that, for some $\mathfrak{m},$ $\mathbf{j}$
 \be\label{i-l}\bigwedge_{n\in \N} n|\mathfrak{m}\ \&\ \mathfrak{m}^2=\lfr\ \&\ \mathbf{m}|\mathbf{j}\ \&\ \mathfrak{j}^2=\ii\ \&\  \ii|(\pp-1)\ee

  and 
 \be\label{ep}\epsilon^\eta=1\,\mbox{mod}\, \pp\mbox{ iff }(\pp-1)|\eta\ee

Notation
$$\acl(X)= \mbox{ the algebraic closure of subset $X\subset \F_\pp$ in } \F_\pp.$$

$$\ZZ[N:M]:=\{ z\in \ZZ: \ N\le z\le M\}$$ 
\epk
\bpk \label{i^2}{\bf Lemma.} {\em We may assume in addition to (\ref{i-l}) that }
 $$  \ii^2\neq \pp-1 \Rightarrow \bigwedge_{a_1,\ldots,a_k\in \Z}\
\ \ii^k+ a_1\ii^{k-1}+\ldots+a_k\neq 0\,\mathrm{mod}\,\pp  
 $$
{\bf Proof.} The set of conditions on $\ii$ on the right hand side of $\Rightarrow$ is countable as is the set of conditions in (\ref{i-l}). Thus these can be realised together with $\mathfrak{j}$ and $\lfr$ in the $\omega$-saturated structure $\ZZ.$ $\Box$
  
\epk
\bpk \label{alpha} {\bf Lemma}. {\em Asssuming $\pp,\ii$ satisfying (\ref{i-l}) and \ref{i^2} are fixed, 
there is an $\epsilon\in \F_\pp$ which along with (\ref{ep}) satisfies
\be \label{ep+} \epsilon^\frac{\pp-1}{\ii}\notin \acl(\Q[\ii]) \ee
that is  
$$\alpha:=\epsilon^\frac{\pp-1}{\ii}=\exp_\pp(\frac{\pp-1}{\ii})$$ is transcendental in $\F_\pp$ over $\ii.$}

{\bf Proof.}  
 Indeed, $\epsilon^\frac{\pp-1}{\ii}$ is an element of (non-standard) order $\ii$ in the multiplicative group $\F_\pp^\times$ of order $\pp-1.$ The set of integers $r\in \ZZ[0:\ii]$ co-prime with $\pp-1$ is an infinite definable set and for each such $r$ the element $\epsilon^r$ satisfies (\ref{ep}) in place of $\epsilon.$ At the same time  $(\epsilon^r)^\frac{\pp-1}{\ii}$ takes different values for distinct $r.$ Since $\ZZ$ is $\omega$-saturated, there is  an $r$ such that
  $(\epsilon^r)^\frac{\pp-1}{\ii}\notin \acl(\Q[\ii]).$ Taking $\epsilon$ to be $\epsilon^r$ proves the claim.
 
\epk
\bpk We write $(\ZZ; \Omega_\pp)$ for the structure of nonstandard arithmetic with extra parameters $\pp,\epsilon, \ii.$
Note that $\Omega_\pp$ also determines $\mfr$ and $\mathbf{j}$ when (\ref{i-l}) is satisfied.

 Call a $k$-tuple $\mathsf{L}\in \ZZ^k$ 
{\bf finite-generic}  with respect to $\Omega_\pp$ (f,-g. for short)
   if for any $\Omega_\pp$-formula $\Phi(\bar{x})$ 
\be\label{Phi} \mbox{For each }\bar{l}\in \N^K,\  (\ZZ;\Omega_\pp)\vDash  \Phi(\bar{l}) \ \ \Rightarrow\ \  (\ZZ;\Omega_\pp)\vDash \Phi(\mathsf{L})\ee

Note that the set of formulas satisfying the antecedent of (\ref{Phi})  is a $k$-type.  We call a formula $\Phi$ as in (\ref{Phi}) an f.-g. formula. 

\medskip

\bpk\label{propert}{\bf Properties.} 

1. $\bar{n}\in \N^k$ then  $\bar{n}$ is f.-g.

2. $(L_1,L_2)$ is f.-g. then $L_1$ and $L_2$ are f.-g.

3. $L\in \ZZ^k$ f.-g. and $f(X)$ an $\Omega_\pp$-definable function such that, $f(\N^k)\subseteq \N.$ Then $f(L)$ is f.-g. and
$f(L)\ge 0.$ 

4. Let  $\beta$ be a definable elemenent in $(\ZZ; \Omega_\pp).$
If $x=L$ is a solution of non-trivial equation 
$p(x,\beta)=0\mod \pp,$  for $p(x,y)\in \Z[x,y],$  then
$L$ is not f.-g.

{\bf Proof.} 1. and 2. immediate by definition. 

3.  Let $\Phi(x)$ be an f.-g. formula. Then $\Phi(f(\bar{y}))$ is f.-g. since \linebreak
$(\ZZ;\Omega_\pp)\vDash \Phi(f(\bar{n}))$ for all $\bar{n}\in \N^k.$
Thus  $\vDash\Phi(f(L)).$ So $f(L)$ is f.-g.. Note that $x\ge 0$ is a f.-g. formula, so $f(L)\ge 0.$

4. There are only finitely many solutions of 
$p(x,\beta)=0\mod \pp\ \& \ 0\le x<\pp,$ so $L$ not f.-g.. 
$\Box$    

\epk

\epk
\bpk {\bf Remark.} We would need to consider also the structure on the field $\C$ with a predicate distinguishing $\Z$ and constants symbols $\ii,$ $\epsilon$ and $\pp,$ call it $(\C, \Z, \Omega_\pp),$ together with its non-standard version $({^*\C}, \ZZ; \Omega_\pp).$ It is easy to check that

{\em  A subset $P\subseteq \ZZ^n$ is definable in $({^*\C}, \ZZ; \Omega_\pp)$ if and only if it  is definable in $(\ZZ; \Omega_\pp).$}

\medskip

Thus in the definition above we can consider either of the structures.
   
\epk
Let $\mathcal{F}$ be the set of all $\Omega_\pp$ definable functions
$f: \ZZ\to \ZZ$ such that $f(\N)\subseteq \N.$

\medskip

{\bf Remark.} We don't know if $\mathcal{F}$ is a proper extension of $\mathcal{F}_0,$ the set of all functions 0-definable in arithmetic. 

However, in more general setting, there is a nonstandard $\mathfrak{q}\in \ZZ$ such that 
the set $\mathcal{F}_\mathfrak{q}$ of all functions $g: \ZZ^k\to \ZZ$ defined in $(\ZZ;+,\cdot, \mathfrak{q})$ and satisfying the condition $g(\N^k)\subset \N$ is bigger than  $\mathcal{F}_0.$ Indeed, set 
$$g(n):=\left\lbrace\begin{array}{ll}
p_n, \mbox{ if } p_n|\mathfrak{q}\mbox{ ($n$-th prime)} \\
0,\mbox{ \ otherwise}
\end{array}\right. $$
Using saturatedness one can, for any set $Q$ of standard primes $p_n\in \N,$ find $\mathfrak{q}\in \ZZ$ such that, for all $n\in \N,$ 
$$ p_n\in Q  \Leftrightarrow p_n|\mathfrak{q}$$ 
Clearly, there are $Q\subset \N$ which are not 0-definable in arithmetic. So exists $\mathfrak{q}$ and $g\in \mathcal{F}_\mathfrak{q}\setminus \mathcal{F}_0.$

\bpk\label{PropM}{\bf Proposition.}
{\em  Assume (\ref{i-l})-(\ref{ep+})  for some $\Omega_\pp$ and  $\lfr.$  Assume $\ii^2+1\neq \pp.$

Then there exist a finite-generic  
$\lfr$ satisfying  (\ref{i-l}) 
together with the condition: for each $f\in \mathcal{F}$
 \be \label{i^2+}  \forall a_0\ldots a_k\in \ZZ[-f(\lfr): f(\lfr)]\
a_0\neq 0\to \ a_0\ii^k+ a_1\ii^{k-1}+\ldots+a_k\neq 0\,\mathrm{mod}\,\pp  
 \ee
In particular, independently on the condition $\ii^2+1\neq \pp,$ for each $f\in \mathcal{F}:$
\be \label{ii>l} \ii > f(\lfr)\ \mbox{ and }\ \frac{\pp-1}{\ii}> f(\lfr)
\ee
}

{\bf Proof.} Let, for $M\in \N,$ 
$$P_M(l):= \bigwedge_{n\le M} \exists m, j :\ n|m \ \& \ m^2=l\ \& \ m|j\ \& \ j^2=\ii\ \& \ \ii|(\pp-1).$$
Clearly $\bigwedge_{M\in \N}P_M(l)$ is a type realisation $\lfr$ of which satisfies 
   (\ref{i-l}). One can see that for any standard $n\ge M,$ $$(\ZZ; \Omega_\pp)\vDash P_M((n!)^2)$$ 
   That is the $\Omega_\pp$-formula 
   $$\Psi_M(n):= \ n\ge M \to P_M((n!)^2)$$
  holds for all $n\in \N,$ that is $\Psi_M$ is an f.-g,-formula. Let $\mathfrak{n}\in \ZZ$ be an infinite f.-.g. number.  In particular, $(\ZZ; \Omega_\pp)\vDash \bigwedge_{M\in \N} \Psi_M(\mathfrak{n}).$ Let $\lfr:= (\mathfrak{n}!)^2,$ which is also an infinite f.-.g. number by 
3. of  \ref{propert}.


 Let
$\Phi_f(x)$ be the $\Omega_\pp$-formula stating   (\ref{i^2+}) when $x:=\lfr.$ 
Lemma \ref{i^2} states that this is a f.-g. formula. Then 
 $(\ZZ;\Omega_\pp)\vDash \Phi_f(\lfr)$ by our definition of $\lfr,$ which proves (\ref{i^2+}). 
 
 In case $\ii^2+1\neq\pp$ (\ref{ii>l}) follows immediately. Otherwise we use the fact that $\ii>f(n)\ \&\ \frac{\pp-1}{\ii}>f(n)$ for all $n\in \N.$ 
 $\Box$

\epk
\bpk\label{notation} {\bf Notation/Corollary.} We choose $\lfr,$ $\epsilon$ and $\ii$ so that  (\ref{i-l})-(\ref{ep+}) and (\ref{i^2+}) hold and $\lfr$ is finite-generic.

 
Denote
$$\Oo(\mathcal{F}):= \bigcup_{f\in \mathcal{F}}  \ZZ[-f(\lfr): f(\lfr)],$$
$$\Oo(\lfr):= \bigcup_{m\in \N}  \ZZ[-m\lfr: m\lfr],$$
$$\uuu:=\frac{\pp-1}{\ii\lfr}\mbox{ and }\vvv:=\frac{\pp-1}{\lfr}=\ii \uuu.$$
\epk

\bpk\label{scales} {\bf Corollary.} {\em $\Oo(\mathcal{F})$ is a convex subring of $\ZZ$ containing $\Oo(\lfr)$ and closed under every
$f\in \mathcal{F}.$ In particular,
$$\Oo(\mathcal{F})\prec \ZZ$$
in the language of rings.

\medskip

Also
$$\Oo'(\mathcal{F})\ <\ \Oo'(\mathcal{F})\cdot \uuu\ <\ \Oo'(\mathcal{F})\cdot \vvv$$ 
for $\Oo'(\mathcal{F}):=\Oo(\mathcal{F})\setminus \{ 0\}.$

}

{\bf Proof.} The first follows from the definition and the fact that $\Oo(\mathcal{F})$ is closed under all 0-definable maps in the structure $(\ZZ; +,\cdot)$ with definable Skolem functions.

The second is a corollary of  (\ref{ii>l}). $\Box$

\medskip

\epk
\bpk {\bf Definition} 
$$\Sss:= \exp_\pp\{\Oo(\lfr)\cdot\vvv\}$$
$$\RR_+:=\exp_\pp\{\Oo(\lfr)\cdot\uuu\}$$
 
Note that since $\exp_\pp((n+\lfr)\cdot \vvv)=\exp_\pp(n\cdot\vvv)$
$$\Sss:= \exp_\pp\{\Oo(\lfr)\cdot\vvv\}$$ 
\epk

\section{Embedding 
 into ${^*\C}$.}

\bpk \label{3.1}
Let $\pi'\in \CC$ stand for a real number, possibly non-standard, such that $\e^{\pi'}$ is transcendental (assuming Schanuel's conjecture we can take $\pi':=\pi$). 

\medskip

We consider $\U=\U_{\pp,\lfr}$ as an $\Oo(\lfr)$-module. Note that $\ZZ\subset \CC$ and thus $\Oo(\mathcal{F})\subset \CC$ and 
acts on $\CC$ by multiplication.

  First we define the map
$${\mathrm{I}}_\U: \U\to \CC$$
 on the 2-elements set $\{  \uuu, \vvv\}\subset \U$ (see \ref{notation}):
 $${\mathrm{I}}_\U:\ \ \uuu\mapsto -\frac{2\pi'}{\lfr}; \ \ \vvv\mapsto -\frac{2\pi \mathrm{i}}{\lfr}$$
where $\mathrm{i}=\sqrt{-1}.$
 
For $\alpha,\beta\in \Oo(\lfr)$ define
 \be\label{IU} {\mathrm{I}}_\U:\ (\alpha\cdot \uuu+\beta\cdot \vvv)\mapsto -\frac{1}{\lfr}(2\pi'\alpha + 2\pi \mathrm{i}\beta).\ee 
 Note that by \ref{PropM}, $\uuu$ and  $\vvv$ are linearly independent over  $\Oo(\lfr)).$ 
 Hence the map is well-defined and invertible.
 
 Let $\U(\lfr)$ be the 2-dimensional $\Oo(\lfr)$-submodule of $\U$ spanned by $\{ \uuu,\vvv\}$ and let $\CC(\lfr)\subset \CC$ be the  $\Oo(\lfr)$-submodule spanned by $\{ \frac{2\pi'}{\lfr}, \frac{2\pi \mathrm{i}}{\lfr}\}$

(\ref{IU}) defines   
$${\mathrm{I}}_\U:\ \U(\lfr) \to \CC(\lfr)$$ 
as an $\Oo(\lfr)$-linear isomorphism.

\epk
\bpk\label{defIF} Consider elements
$$1,\ \exp_\pp(\uuu),\  \exp_\pp(\vvv)\in \F_\pp$$
($1$ is $1_{\mathrm{mod}\,\pp}$ of $\F_\pp$)
and define a partial map $\mathrm{I}_\F: \F_\pp\to \CC$ on the three points:
$$\mathrm{I}_\F: \ 1\mapsto 1,\  \exp_\pp(\uuu)\mapsto \e^{-\frac{2\pi'}{\lfr}},\ \ \exp_\pp(\vvv)\mapsto \e^{-\frac{2\pi \mathrm{i}}{\lfr}}$$
and further, for any $a\in \Oo(\mathcal{F})$ and $s,r\in \Oo(\lfr),$
$$\mathrm{I}_\F: \ a\cdot 1\mapsto a\cdot 1,\  \exp_\pp(r\uuu)\mapsto \e^{-\frac{2r\pi'}{\lfr}},\ \ \exp_\pp(s\vvv)\mapsto \e^{-\frac{2s\pi \mathrm{i}}{\lfr}}.$$

This is internally definable over respective elements   by our assumptions. Moreover, if $a=a_i,  s=s_i$ and $r=r_i$ represent elements of internally definable sequences, $i\in I\subset \ZZ[0:\lfr^m]$ then one can definably extend to sequences
$$\mathrm{I}_\F: \ a_i\cdot 1\mapsto a_i\cdot 1,\   \exp_\pp(r_i\uuu)\mapsto \e^{-\frac{2r_i \pi'}{\lfr}},\ \ \exp_\pp(s_i\vvv)\mapsto \e^{-\frac{2s_i \pi \mathrm{i}}{\lfr}}.$$
 
Let $$\Oo(\F):=\{ \sum_{i\in I} a_i   \exp_\pp(r_i\uuu)\exp_\pp(s_i\vvv): \ I\subseteq \ZZ[0:\lfr^m] \mbox { definable}\}$$
 $$\Oo(\CC):=\{ \sum_{i\in I} a_i   \e^\frac{-2r_i \pi'}{\lfr}\e^\frac{-2s_i \pi \mathrm{i}}{\lfr}: \ I\subseteq \ZZ[0:\lfr^m] \mbox { definable}\}$$
These are rings containing $\Oo(\mathcal{F})$ and closed under internally definable summation. 
 
Extend $$ \mathrm{I}_\F: \Oo(\F)\to \Oo(\CC)  $$ accordingly.

\epk
\bpk\label{IF} {\bf Lemma.} {\em $\mathrm{I}_\F$ is well-defined and bijective.}

{\bf Proof}. It is enough to prove that 
\be\label{Ph} \sum_{i\in I} a_i   \exp_\pp(r_i\uuu)\exp_\pp(s_i\vvv)=0 \mbox{ iff }
\sum_{i\in I} a_i   \e^\frac{-2r_i \pi'}{\lfr}\e^\frac{-2s_i \pi \mathrm{i}}{\lfr}=0.\ee
Note that the left hand side equality  
 can be expressed by a formula
$\Psi(\lfr,I)$ in $(\ZZ;\Omega_\pp)$. 

Similarly, the right   hand side equality  
 can be expressed by a formula
$\dot{\Psi}(\lfr,I)$ in  
$({^*\C}, \ZZ; \Omega_\pp).$

Consider the formula
$$\Phi(\lfr):= \forall I\subset \ZZ[0:\lfr^m]\ {\Psi}(\lfr,I) \leftrightarrow \dot{\Psi}(\lfr,I).$$
We claim that, for all $n\in \N,$

 $$({^*\C};\ZZ, \Omega_\pp)\vDash \Phi(n).$$
 Indeed, the internally definable sequences on $\Z[0:n^m]$ have finite values of parameters $a_i, r_i, s_i$ and $\mfr$
and thus $\Phi(n)$ expresses the fact that the algebraic dependence 
of some $\e^\frac{2r_i \pi'}{k},$ $\e^\frac{2s_i \pi \mathrm{i}}{k}$ with coefficients $a_i$ 
takes place if and only if  the algebraic dependence 
 with coefficients $a_i$ of 
$\exp_\pp(\frac{r_i(\pp-1)}{k\ii}),$ $\exp_\pp(\frac{s_i(\pp-1)}{k})$
takes place. And note that both $\e^\frac{2s_i \pi \mathrm{i}}{k}$ and 
$\exp_\pp(\frac{s_i(\pp-1)}{k})$ are standard roots of unity in respective fields while $\e^\frac{2r_i \pi'}{k},$ 
$\exp_\pp(\frac{r_i(\pp-1)}{k\ii}),$ are transcendental elements. It follows that the equalities on the both side of (\ref{Ph}) can only happen when all $a_i=0.$ Which proves $\Phi(n)$ and hence $\Phi(\lfr).$ $\Box$

\epk
\bpk \label{aclO} {\bf Remark.} $$\acl(\Oo(\CC))\neq \CC.$$
Indeed,  elements  of  $\Oo(\CC))$ are definable over $\lfr, \e^\frac{\pi'}{\lfr}$ and $\e^\frac{\pi i}{\lfr}$ in $(\CC; \ZZ,\Omega_\pp)$ and thus can be reduced  to elements of the ring $\Oo_l\subset \C$ generated  by roots of unity of order $l$ and elements of the form  $\e^\frac{2\pi' i}{l}$ over $l\in \N$ in a structure of the form $(\C; \Z,\Omega_p)$.
Hence any element of $ \acl(\Oo(\CC))$ is in the ultraproduct of  $\acl (\Oo_l)$ such that $l$ are the restrictions of $\lfr$ to the index set of the ultraproduct. The statement follows from the fact that  $\acl( \Oo_l)\neq \C.$ $\Box$
\epk
The Lemma together with the statement in \ref{3.1} prove:

\bpk\label{Main1} {\bf Theorem.} {\em The maps
$${\mathrm{I}}_\U:\ \U(\lfr) \to \CC(\lfr)$$
and
$$ \mathrm{I}_\F: \Oo(\F)\to \Oo(\CC)  $$
are bijections preserving internally definable summation and commute with respective exponentiation maps
$$ \exp_\pp: \U(\lfr)\to \Oo(\F)\mbox{ and }\exp: \CC(\lfr)\to \Oo(\CC)  $$

}

\epk
\bpk Define $\F\subset \F_\pp$ to be the fraction field of $\Oo(\F).$

\medskip

Theorem \ref{Main1} implies that $\mathrm{I}_\F$ extends to the embedding
\be \label{IF1}\mathrm{I}_\F:\ \F \hookrightarrow \CC\ee
\epk

 \bpk {\bf Corollary.}
 {\em There is an internally definable notion of {\bf complex conjugation} 
 $z\mapsto \bar{z},$  $\bar{z}=z\inv,$  for $z\in {'\Ss'},$ and $\bar{y}=y$ for $y\in \RR_+$  and this determines an automorphism $x\mapsto \bar{x}$ on $\F.$ 
 
 Moreover, $$\mathrm{I}_\F(\bar{x})=\overline{\mathrm{I}_\F({x})}$$
 (complex conjugation in $^*\C$ on the right).}
 
 
{\bf Proof.} $\mathrm{I}_\F$ maps $\Sss$ to $\Ss$ and $\RR_+$ to positive reals of $\CC,$ thus complex conjugation is correctly defined  on $\Sss$ and $\RR_+.$ The extension to $\F$ is by internally definable summation and hence reduces readily to finite sums, which satisfies the algebraic identities of complex conjugation. $\Box$

\epk

Set $$\lm_\F:= \mathrm{st}\circ \mathrm{I}_\F$$

\bpk\label{limF} {\bf Lemma.}
$$\lm_\F: \ \F \twoheadrightarrow \bar{\C}$$

{\bf Proof.} The surjectivity of $\lm_\F$ follows from the facts:
$$\lm_\F: \{ \exp_\pp(r\uuu): \ r\in \Oo(\lfr)\}=\R_+$$ 
$$\lm_\F: \{ \exp_\pp(r\vvv): \ r\in \Oo(\lfr)\}=\Ss$$
and
$$\lm_\F(\lfr)=\infty.$$
$\Box$
\epk

\bpk \label{iianda} {\bf Lemma.} Let $\hat{\ii}=\ii_{\mbox{mod} \pp}\in \F_\pp.$ Then

$$\hat{\ii}^2=-1\mbox{ or } \hat{\ii}\notin \acl(\F).$$

{\bf Proof.} Suppose towards a contradiction that $\hat{\ii}\in \acl(\F).$
It implies that 
\be \label{cc} c_0\hat{\ii}^k+c_1\hat{\ii}^{k-1}+\ldots +c_k=0\ee
for some $c_0, c_1,\ldots c_k\in \Oo(\F),$ $c_0\neq 0,$ that is of the form \linebreak $ \sum_{i\in I} a_i   \exp_\pp(r_i\uuu)\exp_\pp(s_i\vvv)$ each, i.e.  $c_i=c_i(\lfr)$  internally definable in $(\ZZ;\Omega_\pp)$ over $\lfr.$ Note that if we substitute $n\in \N$ in place of $\lfr,$ $c(n)$ is in the ring generated by $\exp_\pp(\frac{r_i(\pp-1)}{k\ii}),$ $\exp_\pp(\frac{s_i(\pp-1)}{k})$ with finite $r_i,s_i$ and $k$. That is $c(n)$ is algebraically dependent on $\alpha=\exp_\pp(\frac{\pp-1}{\ii}).$ It follows  that, 
 \be\label{negcc} \mbox{for all } n\in \N:\ \ c_0(n)\hat{\ii}^k+c_1(n)\hat{\ii}^{k-1}+\ldots +c_k(n)\neq 0\ \vee \ \bigwedge_{0\le i\le k}c_i(n)=0\ee
for otherwise $\hat{\ii}\in \acl(\alpha)$ in contradiction with \ref{alpha} and \ref{PropM}.

Since $\lfr$ is finite-generic  (\ref{negcc}) implies the negation of (\ref{cc}), the contradiction which proves our statement.$\Box$ 
\epk
Recall that $$\mathrm{i}=\e^{\frac{\pi i}{2}}.$$
\bpk\label{3.8} {\bf Corollary.} {\em For some $\mathrm{i}'\in \CC$ such that $\mathrm{i}'-\mathrm{i}$ is an infinitesimal in $\CC,$
the embedding $\mathrm{I}_\F$ of (\ref{IF1}) extends to the embedding $$\mathrm{I}_\F:\F_\pp\to {\CC};\ \mbox{ so that } \hat{\ii} \mapsto \mathrm{i}'$$

Thus, 
\be\label{limF+} \lm_\F: \ \F_\pp \twoheadrightarrow \bar{\C}\ee
 \be\label{lmj} \lm_\F: \hat{\ii}\mapsto \e^{\frac{\pi i}{2}}\mbox{ and } \hat{\jj}\mapsto \e^{\frac{\pi i}{4}}.\ee
}

{\bf Proof.} If ${\hat{\ii}}^2=-1$ then set $\mathrm{i}'=\mathrm{i}.$
Otherwise, pick $\mathrm{i}'\in \CC$ in the infinitesimal neighborhood of $\mathrm{i}$ but not in $\mathrm{I}_\F(\acl(\F)),$ which exists because of \ref{aclO}.

The extension of  $\mathrm{I}_\F$ to an embedding of $\F_\pp$ is by the routine algebraic construction using the fact that $$\mathrm{tr.deg}_\F \F_\pp\le  \mathrm{tr.deg}_\F \CC$$
and $\CC$ is algebraically closed. $\Box$

\medskip

{\bf Remark.} Note that in terms of  the embedding $\ZZ[0: \pp-1]\hookrightarrow \F_\pp$ one may identify 
 $\hat{\ii}=\ii,\ \hat{\jj}=\jj$ 
 and $$\lm_\F: {\ii}\mapsto \e^{\frac{\pi i}{2}}\mbox{ and } {\jj}\mapsto \e^{\frac{\pi i}{4}}.$$

\medskip

{\bf Remark} An immediate consequence of the properties of the standard part map is that
$\lm_\F$ as defined in (\ref{limF+}) is
 a place of fields, that is there is 
a local ring $\F^0$ of $\F$ such that the restriction $$\lm_\F: \F^0\twoheadrightarrow \C$$   is a homomorphism of rings and, for $x\in \F\setminus \F^0,$ $\lm_\F(x)=\infty.$ 
\epk
\bpk {\bf Theorem.} 
{\em There is an additive surjective homomorphism
$$\lm_\U:  \U(\lfr)\twoheadrightarrow \C$$
such that, for $u\in \U(\lfr)$ and $x\in \F$
$$\lm_\U (u)=\mathrm{st}\circ  \mathrm{I}_\U(u)$$
$$\lm_\F (x)=\mathrm{st}\circ  \mathrm{I}_\F(x)$$
$$\exp(\lm_\U(u))= \lm_\F(\exp_\pp(u))$$
and $$\lm_\U( \ii u)=\mathrm{i}\, \lm_\U(u).$$
}

{\bf Proof.} Let
$$\lm_\U (u):=\mathrm{st}\circ  \mathrm{I}_\U(u)$$
which is well-defined for $u\in \U(\lfr)$ by \ref{Main1}, and satisfies 
the required commutativity condition by the same Theorem. 

(\ref{IU}) implies $\lm_\U( \ii u)=\mathrm{i}\, \lm_\U(u)$

The surjectivity of $\lm_\U$ follows from the fact that $$\mathrm{st}(\CC(\lfr))=\C$$ since by construction 
$$\mathrm{st} (\frac{2\pi'\cdot \Oo(\lfr)}{\lfr})=\R\mbox{ and }  \mathrm{st} (\frac{2\pi i \cdot \Oo(\lfr)}{\lfr})=\imath \R.$$  
 $\Box$
 
\epk
\bpk\label{distance} {\bf Order, distance and continuity}. Following \ref{Main1} we are going to assume
\be\label{FtoC} \F\subset {^*\C}.\ee
This allows to consider  the inequality $\le$ on the reals of $\F,$ the restriction of the internally definable relation $\le$ on the reals of  ${^*\C}.$


\medskip

More generally, suppose $X\subset \U^M$ is a definable set with the structure of a $\QQ_+$-valued length-metric, namely
there are internally definable ternary predicates on $\QQ_+\times X^2,$ written as $d_q(x,y),$ 
which are interpreted as ``the distance between $x$ and $y$ is $\le q$''.  Since  $X$ is pseudo-finite the {\bf distance} $\mathrm{dist}(x,y),$  equal to the minimum length of the paths between the two points,  is a well-defined value in $\QQ_+$.

\medskip 

A map $g: X\to \F$ will be called (Lipschitz) {\bf pseudo-continuous} (with derivative bounded by $c$) if
there exists positive $c\in \Q$ such that for any $x_1,x_2\in X,$ 
$$d_{\frac{1}{\lfr}}(x_1,x_2)\to |g(x_1)-g(x_2)|\le \frac{c}{\lfr}.$$

\epk
\bpk\label{bound} {\bf Lemma.} {\em Let $$g: X\to \F$$
be  pseudo-continuous with derivative bounded by $c.$ Then for all $M\in \Oo(\mathcal{F})$ for all $z_1,z_2\in X$}
$$ d_{\frac{M}{\lfr}}(z_1,z_2) \ \to\   |g(z_1)-g(z_2)|\le c\frac{M}{\lfr}.$$

{\bf Proof.} Immediate from definition by induction on $M.$ $\Box$
\epk

 \section{States and the Hilbert space}

\bpk \label{Dom}
We will assume that $\Dom\subseteq \U^M$ introduced in \ref{LE} is a set  with some family $\Omega_\Dom$ of internally definable  relations on it. In particular, $\Omega_\Dom$ contains the predicates for the {\bf $\QQ$-metric structure} on $\Dom:$
this is a family of binary predicates $d_q(x,y),$ $0<q\in \QQ,$ with intended interpretation ``the distance between $x$ and $y$ is $\le q$''.


 \medskip

A  {\bf  state on} $\Dom$ is an internally definable 
map
$$\varphi: \Dom\to \F.$$

A {\bf ket}-state on $\Dom$ is a pseudo-continuous state  on $\Dom$.

\medskip

Set $\HH^-_\Dom(\F)$ to be the multisorted structure each sort being a state on
 $\Dom.$  
Equivalently, $\HH^-_\Dom(\F)$ is the set of all states on $\Dom.$

\medskip

We denote $\NN:=\#\Dom$, the (non-standard) number of elements in $\Dom.$  

\medskip 
 
We say that $\Dom$ is {\bf tame} if there is an internally definable embedding
$$\Dom\to \ZZ[0:\lfr^m].$$  
It follows
$$\NN\le \lfr^m$$
for tame $\Dom.$

\epk

\bpk \label{positions}
Define now a special kind of states, the position states 
$$\uu[r]:\ \Dom\to \{ 0,1\}\subset {\Sss}; \ \ r\in \Dom$$
$$\uu[r](x)=\delta(r-x):=\left\lbrace\begin{array}{ll}
1, \mbox{ if }r=x\\
0,\mbox{ otherwise}
\end{array}\right.$$
Define, for a ket-state $\varphi$ the inner product with $\uu[r]$
$$\la \varphi | \uu[r]\ra:=\varphi(r)$$
and also 
$$\la \uu[x] | \uu[r]\ra:= \delta(r-x)$$
 
 It is immediate from the definition that
 
\epk
\bpk
{\bf Lemma.} {\em Assume $\Dom$ is tame. Then $\HH^-_\Dom$ can be given the structure of $\F$-linear space with an inner product defined as
\be \label{<>}\la \varphi| \psi\ra:=\sum_{r\in \Dom} \varphi(r)\cdot \bar{\psi}(r)\in \F\ee

The definition is consistent with $\varphi$ or/and $\psi$ being position states.

The product is Hermitian, that is satisfies the sesquilinearity condition and is positive definite.}


{\bf Proof.} By assumptions $r\mapsto \varphi(r)\cdot \bar{\psi}(r)$ can be identified as an internally definable sequence from a subset of $\ZZ[0:\lfr^m].$ Thus the sum (\ref{<>}) is well-defined and belongs to $\F.$
The rest is immediate by definitions. $\Box$ 

\medskip

Define the square of the {\bf norm}
$$|\varphi|^2:={\la \varphi|\varphi\ra}\in \R_+(\F).$$

In terms of $\F$-logical values the inner product estimates the equality
``$\varphi=\psi$'' for $\varphi,\psi$ of norm 1. Indeed, ``$\psi=\psi$'' is given value $1.$
\medskip

 Set {\bf the range} of $\varphi$ to be:
$$\mathrm{Range}(\varphi):=\{ r\in \Dom:\ \la \varphi| \uu[r]\ra \neq 0\}$$
\epk

\bpk {\bf Lemma.} {\em Suppose
  $|\varphi(r)|\le \eta$ for all $r\in \Dom.$ Then
\be\label{range} |\varphi|^2:={\la \varphi|\varphi\ra}\le \eta\cdot \#\mathrm{Range}(\varphi)\ee
}

{\bf Proof.} By definition $$\la \varphi|\varphi\ra=\sum_{r\in \Dom} |\varphi(r)|^2=\sum_{r\in \Dom} |\la \varphi(r)| \uu[r]\ra|^2$$
and since $\Dom$ is tame and the sum is internally definable in $^*\C$ we can lift the required inequality from the finite to pseudo-finite summation. $\Box$

\
\epk
\bpk Let   $\{ \psi_i: s\in I\}$ be an internally definable family of states on $\Dom$ over an internally definable set $I$.

Then the sum
$$S=\sum_{i\in I} \psi_i$$
is internally definable,
the value $S$
of the sum is an element  of the ultraproduct $\prod_{p\in I} \F_p/D$ such that $S(p)=\sum_{i\in I(p)} \psi_i$ (finite sum) in $\F_p$ along the ultrafilter.

Define $\HH_\Dom(\F)$ to be the smallest $\F$-linear subspace of
$\HH^-_\Dom(\F)$  closed under taking  definable sums.

\epk
\bpk  
It follows from definitions that the set of all  position states 
of $\Dom$  
$$\LL_\uu:=\{ \uu[r]:  \ r\in \Dom\}$$
is  definable.

\medskip

{\bf Lemma.} {\em The set of position states forms a basis of $\HH_\Dom$ with regards to definable summation. This basis is
  orthonormal.
$$\dim \HH_\Dom=\NN$$  }
  
  {\bf Proof.} For an arbitrary $\psi\in \HH_\Dom,$
  $$\psi=\sum_{r\in \Dom} \psi(r)\cdot \uu[r].$$
  
   
  \epk
  \bpk\label{inn} {\bf Corollary.} {\em The definition of inner product  in (\ref{<>}) is applicable to any pair of states in $\HH_\Dom(\F)$ with tame $\Dom.$}
  
  \medskip
  
 From now on we consider  $\HH_\Dom:=\HH_\Dom(\F)$ as an $\NN$-dimensional inner vector space.
  \epk 
 
\bpk\label{operators} {\bf Linear unitary operators on $\HH_\Dom.$}

We are going to consider linear operators 
 $$\A: \HH_{\Dom_1}\to \HH_{\Dom_2}$$
 for $\Dom_1,\Dom_2\subset \U^M$ definable domains.
 
 Call such an operator  definable if
$$\{ \A\uu[r]: \ r\in \Dom_1\}$$
is a definable family  of states on $\Dom_2.$

By definition \be\label{Aa} \A \uu[r]=\sum_{s\in \Dom_2} a(r,s)\cdot \uu[s]\ee
with unitarity condition
\be\label{Aunitary} \sum_{s\in \Dom_2} |a(r,s)|^2=1
\ee
and we also require that the map
$$s\mapsto a(r,s);\ \ \Dom_2 \to \F$$
be pseudo-continuous.

\epk

\bpk \label{product}
{\bf Proposition.} {\em Let $\LL_\psi:=\{ \psi_r: r\in \Dom_1\}$ be a definable basis of $\HH_{\Dom_1}.$ Then 

- there is a definable operator $\A$ such that
$$\A\,\uu[r]=\psi_r; \ r\in \Dom_1$$

- for any definable operator $B$ on $\HH_{\Dom_2}$ the set
$$\B\,\LL_\psi:=\{ \B\,\psi_r: \ r\in \Dom_1\}$$
is definable.

Assuming that $\B$ is unitary the product operator $\B\cdot \A$ is unitary. }

{\bf Proof.} The statement can be equivalently reformulated as a property of finite dimensional vector spaces $\HH(p),$ $p\in P,$ along the ultrafilter $D.$ In this form it is obvious and 
$$
\B\, \psi_s=(\B\cdot \A)\,\uu[s]$$
the product of matrices.
 $\Box$

\medskip

We only consider definable linear operators on $\HH_\Dom.$ In applications (as in \ref{QM}) we do not assume that the operators are invertible but most of the time we deal with operators such as the image of $\LL_\uu$ is of size $\frac{\NN}{k}$ for some finite $k.$ Moreover we have to deal with operators whose domain is a proper subset of $\LL_\uu$ of size $\frac{\NN}{k}$ for some finite $k.$ 
\epk 
\bpk\label{remarkDef} {\bf Remarks.}

1. Looking at states $\psi_r$ and definable set of states $\LL_\psi:=\{ \psi_r: r\in \Dom\}$ as structures, an operator $$\A: \uu[r]\mapsto \psi_r$$   can be seen as an {\bf interpretation} of $\LL_\psi$ in $\LL_\uu$.

 Sums of states  $\varphi+\psi$ are just objects interpretable in $\LL_\varphi\dot{\cup}\LL_\psi$ (with the intended logical meaning ``$\varphi$ or $\psi$'').

\medskip

2. In regards to model-theoretic formalism,
$\HH_\Dom$ is not considered to be a universe of a structure. Instead, we consider the multisorted structure on sorts $\LL_\psi$ with linear maps between these, together with the $\F_\pp$-linear space $\HH_\Dom$ interpretable in the sorts.

\medskip

3. The multisorted structure $\HH_\Dom$ is by construction interpretable in $(\U, \F_\pp; \Omega_\pp).$ 
 
\epk
\bpk {\bf Proposition.} (a) {\em A bijective transformation
$\sigma: \Dom\to \Dom$ induces the  linear unitary  transformation of $\HH_\Dom:$
$$U_\sigma: \psi\mapsto \psi^\sigma; \ \ \psi^\sigma(r)=\psi(\sigma(r))$$ 
}

(b) {\em Let 
$$\G(\Dom):=\mathrm{Aut}(\Dom)\mbox{ and }\ \G(\HH_\Dom):= \{ U_\sigma: \sigma\in \G(\Dom)\} $$

Then $ \G(\HH_\Dom)\subseteq\mathrm{SU}_\NN(F)\subset \mathrm{GL}_{\NN}(\F)$  is the unitary linear group and $$ \sigma \mapsto U_\sigma$$
is an injection.

In other words, the structure on $\Dom$ is reflected in the algebra of linear operators on $\HH_\Dom.$
}.

{\bf Proof.} (a) Linearity: for $a_1,a_2\in \F,$ $\psi_1,\psi_2\in \HH_\Dom,$
$$(a_1\cdot \psi_1+a_2 \cdot \psi_2)^\sigma(r)=(a_1\cdot \psi_1+a_2 \cdot \psi_2)(\sigma(r))=(a_1\cdot \psi_1^\sigma(r)+a_2 \cdot \psi_2^\sigma(r).$$

Unitarity: it is enough to prove it for basis $\uu[r]: r\in \Dom.$
$$\la \uu[r]|\uu[s]\ra =\delta_{r,s}=\delta_{\sigma(r),\sigma(s)}=
\la \uu^\sigma[r]|\uu^\sigma[s]\ra.$$

(b) It is clear that $\sigma\mapsto U_\sigma$ is a homomorphism. Suppose that $U_\sigma$ is in the kernel. Then $\uu^\sigma[r]=\uu[r],$ which implies by definition of $\uu$ that $\sigma(r)=r,$ for all $r\in\Dom,$
that is $\sigma=\mathrm{id}_\Dom.$ $\Box$ 

\epk
\bpk {\bf The dual to $\HH_\Dom$ and model-theoretic $\HH_\Dom^\mathrm{eq}.$}

Every $\psi\in \HH_\Dom$ gives rise to an $\F$-linear map

$$L_\psi: \xx\mapsto \la \xx|\psi\ra, \ \ \HH_\Dom\to \F$$

By definition, $L_\psi$ is uniquely determined by its values on a basis, that is
$$L_\psi^\uu: \uu[r]\mapsto \la \uu[r]|\psi\ra,\ \ r\in \Dom$$
determines the linear map $L_\psi.$ Clearly, $L_\psi^\uu$ is definable and thus we can treat $L_\psi$ (otherwise given by the  Dirac delta-function) as a definable, or interpretable in $\HH_\Dom.$ 

Set $\HH_\Dom^*$ to be the $\F$-vector space of interpretable $\F$-linear maps with a naturally induced Hermitian inner product structure.   
By construction there is a natural embedding
$$\HH_\Dom^*\subset \HH_\Dom^\mathrm{eq};\ \ \psi\mapsto L_\psi$$
where the imaginary sorts are understood in the sense of internal definability.

\epk

\section{Gaussian Hilbert space}
\bpk\label{GHDom} Here we consider states over a domain $\Dom$ which is not assumed to be tame. Until the last subsection  of this section our domain is ``one-dimensional'' that is
$\Dom\subseteq \U.$ In the last subsection we allow domains $\Dom^M\subseteq \U^M$ for $M\in \Oo(\mathcal{F})$ and then the respective Hilbert space is the $M$-th tensor power of the Hilbert space for the present case.
  
We require $\Dom$ to be of the form $$\Dom=\uuu\cdot \ZZ[-\frac{\NN}{2}: \frac{\NN}{2}-1],$$ for some $\uuu\in \U$ and
$\NN\in \ZZ$ such that $\NN| (\pp-1),$ $\sqrt{\NN}\in \ZZ$ and $\mfr| \sqrt{\NN}.$ 

We treat $\ZZ[-\frac{\NN}{2}: \frac{\NN}{2}-1]$ with addition and multiplication as
the residue ring $\ZZ/\NN$ and treat $\Dom$ as a $\ZZ/\NN$-module. 

Let 
$$ \e(\frac{n}{2\NN}):=\exp_\pp((\pp-1)\frac{n}{2\NN}).$$
and
$$\GG(\NN)=\{  \e(\frac{n}{\NN}): \  n\uuu\in \Dom\},$$
the cyclic subgroup of order $2\NN,$ 
\medskip

We consider quadratic forms $f(x,y)=ax^2+2bxy+cy^2,$  $a,b,c\in \Z,$   where $x,y$ run in the ring $\ZZ/\NN$ often represented by   $\ZZ[-\frac{\NN}{2}: \frac{\NN}{2}-1].$

(The assumption that $f$ is over $\Z$ rather than over $\Q$ is not  restrictive because we can treat 
$\e(\frac{f(x,y)}{k\NN})$ as element of $\GG(k\NN).$ )

\medskip

Denote 
$$\e(\frac{f(x,y)}{2\NN}):=\exp_\pp(\frac{\pp-1}{2\NN}f(x,y))\in \GG(\NN),\mbox{ for }x\uuu,y\uuu\in \Dom.$$

Note that  
$$\e(\frac{f(x,y)}{2\NN})=\e(\frac{f(x',y')}{2\NN})\mbox{ if } x=x'\ \&\ y=y' \mod \NN.$$

The following set of {\bf   Gaussian coefficients} play an important role:

$$\Oo\GG(\NN):=\{ \frac{c}{\sqrt{\NN}} \cdot \e(\frac{n}{2\NN}):\ n\uuu\in \Dom,\ kc\in \Oo(\mathcal{F}), \ k\in \Z\}$$

\epk
\bpk {\bf Gauss quadratic sums in $\F_p.$} Let $M\in \N,$ even, $p>2$ and  $\zeta:=\e^\frac{\pi i}{M}\in \C,$ primitive root of unity of order $2M.$
The classical Gauss quadratic sums formula (the basic form)
is
$$\sum_{0<n\le M} \zeta^{n^2}=\sqrt{M}\,\e^\frac{\pi i}{4}.$$
Suppose in addition $M$ is divisible by $4.$ Then $$\zeta^\frac{M}{4}=\e^\frac{\pi i}{4}\mbox{ and }\sum_{0<n\le M} \zeta^{n^2}=\sqrt{M}\,\zeta^\frac{M}{4}$$

Let $\Z[\zeta]\subset \C$ be the ring generated by $\zeta.$ By the above $\sqrt{M}\in  \Z[\zeta].$ Let
$\F_p:=\Z/p$ and let $\xi\in \F_p$ be a primitive $2M$-root of  $1$ in $\F_\pp,$ that is  
$$\xi^{2M}=\hat{1}:=1_{\mathrm{mod} p},\ \mathrm{ord}\,\xi=2M $$
Then there is a ring homomorphism 
$$h:\ \Z[\zeta]\to \F_p\mbox{ such that }\zeta\mapsto \xi$$
Clearly, 
 $$h: \ 1\mapsto \hat{1},\ M\mapsto \hat{M},\ \sqrt{M}\mapsto \sqrt{\hat{M}}$$
 and
 $$\sum_{0<n\le M} \xi^{n^2}=\sqrt{\hat{M}}\cdot\e(\frac{1}{8}),\mbox{ where }
\e(\frac{k}{m}):= \xi^\frac{kM}{m}.$$

Simple algebraic manipulation produce a more general version both for characteristic zero and characteristic $p,$  for $a,b\in \Z,$  assuming $\pm a>0$ and $4a|M,$
 \be\label{GaussSum0}\sum_{0<n\le \frac{M}{|a|}} \e({\frac{an^2+2bn}{2M}})=\left\lbrace\begin{array}{ll}\sqrt{\frac{\hat{M}}{|a|}}\e(\pm\frac{1}{8})\e(-\frac{b^2}{2aM})\mbox{ if } a| b \\
0\mbox{ otherwise }\end{array}\right.\ee 
  (Note that the function $\e(\frac{an^2+2bn}{2M})$ of $n$ has period $\frac{M}{|a|}.$ ) 
 
\epk
\bpk\label{Gsummation}
{\bf Gaussian summation over $\Oo\GG(\NN).$}  


\medskip

In the context of states on $\Dom$ we can write (\ref{GaussSum0}) equivalently, for $\pm A>0,$
\be\label{GaussSum} \sum_{nA\uuu\in \Dom}
\e(\frac{An^2+2Bn}{2\NN})=
\left\lbrace \begin{array}{ll}
\e(\pm \frac{1}{8})\sqrt{\frac{\NN}{|A|}}\cdot  \e(-\frac{B^2}{2A\NN}),\mbox{ if }A|B\\
0,\mbox{ otherwise (e.g. $A=0$)} 
\end{array}\right. \ee
 which can be interpreted as a calculations over $\Dom_\uuu$ in units of scales $A\uuu,$ replacing  $\Dom_\uuu$ by  $A\Dom_\uuu.$ 
 
We say that a  quadratic form over $\Z,$\ $f(u,v)=Au^2+2Buv+Cv^2,$ is {\bf admissible} if  $A,C\le 0$

\medskip

Consider $$f_1(x,y)=A_1x^2+2B_1xy+C_2y^2\mbox{ and } f_2(x,y)=A_2x^2+2B_2xy+C_2y^2$$ be quadratic forms over $\Z,$ $A=A_1+A_2,$ $B=B_1+B_2.$

$$\sum_{nA\uuu\in \Dom}\e(\frac{f_1(n,p_1)}{2\NN})\cdot \e(\frac{f_2(n,p_2)}{2\NN})=$$ $$=\e(\frac{C_1p_1^2+C_2p_2^2}{2\NN})\cdot \sum_{nA\uuu\in \Dom}
\e(\frac{(A_1+A_2)n^2+2(B_1p_1+B_2p_2)n}{2\NN})=$$ $$=\e(\frac{C}{2\NN})\sum_{nA\uuu\in \Dom}
\e(\frac{An^2+2Bn}{2\NN});\ \ C=C_1p_1^2+C_2p_2^2. $$


This can be rewritten with normalising coefficients as
\be\label{sumInner}\sum_{nA\uuu\in \Dom}\frac{1}{\sqrt{\NN}}\e(\frac{f_1(n,p_1)}{2\NN})\cdot \frac{1}{\sqrt{\NN}}\e(\frac{f_2(n,p_2)}{2\NN})= \frac{c}{\sqrt{\NN}}\cdot\e(\frac{C-\frac{B^2}{A}}{2\NN})\in \Oo \GG(\NN)\ee
where 
$$c=\left\lbrace \begin{array}{ll}\e(\pm\frac{1}{8})\sqrt{\frac{1}{|A|}},\mbox{ if }A|B\\
0, \mbox{ otherwise} 
\end{array}\right.$$

\epk

\bpk {\bf Gaussian-ket states on $\Dom$} are  definable sequences of elements of $\F_\pp$ of the form
$$\s_f[p]:= \{ \cc_\s\cdot\e(\frac{f(r,p)}{2\NN}): \ r\uuu\in \Dom\}, \mbox{ where } p\uuu\in \Dom,$$ 
where $f(r,p)$ is an admissible quadratic form over $\Z$
  and  $c_\s\in \Oo\GG(\NN).$
 
 Equivalently these can be written as
symbolic expressions of the form
$$\s_f[p]:=  \cc_\s\,\sum_{r\uuu\in \Dom}\e(\frac{f(r,p)}{2\NN})\,\uu[r],\ \ p\uuu\in \Dom$$
 

\medskip

A Gaussian ket-sort  is
$$\LL_f  :=\{ \s_f[p]: p\uuu\in \Dom\}$$
where  $\s_f$ is as above.

We also consider the positions sort
$$\LL_\uu  :=\{ \uu[r]: r\uuu\in \Dom\}$$

\medskip 

Note that the definition of ket-states allows multiplication  by elements of $\F.$ We consider  $\HH_\Dom$ as the $\F$-linear space finitely generated by elements of Gaussian ket-sorts  $\LL_f $ and position sort  $\LL_\uu.$
\epk
\bpk Now we consider two possible {\bf formal inner product} with values in $\F_\pp$ beween Gaussian ket-states. Let $f_1(x,y)$ and $f_2(x,y)$ be quadratic forms as defined in \ref{Gsummation}.

\medskip

 {\bf formal-Euclidean inner product of ket-states} is defined as
$$\la \s_{f_1}[p_1]\, |\, \s_{f_2}[p_2]\ra_\mathrm{E}=\left\lbrace
\begin{array}{ll}
0,\mbox{ if }A=A_1+A_2=0\\
\cc_1\cc_2\cdot 
\sum_{rA\uuu\in \Dom} \e(\frac{f_1(r,p_1)}{2\NN})\cdot \e(\frac{f_2(r,p_2)}{2\NN}) \mbox{ if }A\neq 0
\end{array}\right.$$  

\medskip

{\bf formal-Hermitian inner product of ket-states} is defined as

$$\la \s_{f_1}[p_1]\, |\, \s_{f_2}[p_2]\ra_\mathrm{E}=\left\lbrace
\begin{array}{ll}
0\mbox{ if }A=A_1-A_2=0\\
\cc_1\cc_2\cdot 
\sum_{rA\uuu\in \Dom} \e(\frac{f_1(r,p_1)}{2\NN})\cdot \e(-\frac{f_2(r,p_2)}{2\NN}) \mbox{ if }A\neq 0
\end{array}\right.$$ 


The Euclidean and Hermitian inner products between a ket-state and a position state are both defined as
$$\la \s_{f}[p]\, |\, \uu[r]\ra:=\cc_\s\cdot \e(\frac{f(r,p)}{2\NN})$$


\epk

\bpk {\bf Lemma.} $$\la \s_{f_1}[p_1]\, |\, \s_{f_2}[p_2]\ra\in  \Oo\GG(\NN)$$
$$\la \s_{f}[p]\, |\, \uu[r]\ra\in \Oo\GG(\NN)$$ both for Hermitian and Euclidean inner product.

{\bf Proof.} Follows directly from (\ref{sumInner}) and calculations above. $\Box$\epk
\bpk Consider a  linear unirary operator $\A$ which is defined on $\LL_\uu$ as
$$\A: \uu[q]\mapsto   \frac{1}{\sqrt{\NN}}\,\sum_{r\uuu\in \Dom}\e(\frac{a(q,r)}{2\NN})\,\uu[r]$$
 $a(u,v)$ a  quadratic form over $\Z$ and $\cc_g\in \Oo\GG(\NN)$
 and acts on $\LL_f$ by pseudo-finite linearity
 \be\label{Ab} \A: \s_f[p]\mapsto \frac{\cc_\s}{\sqrt{\NN}} \sum_{q\uuu\in \Dom} \e(\frac{f(p,q)}{2\NN})\cdot\e(\frac{a(q,r)}{2\NN}) \uu[r]\ee
 
{\bf Remark.} In general  the image $\A\s_f[p]$ may fail to be  pseudo-continuous on $\Dom.$  Namely, the $r$-coordinate  of $\A\s_f[p]$ 
$$\A\s_f[p](r)= \frac{\cc_\s}{\sqrt{\NN}}\sum_{q\uuu\in \Dom} \e(\frac{f(p,q)}{2\NN})\cdot\e(\frac{a(q,r)}{2\NN})$$
is zero outside $k\Dom+d$ for some
 $k=k(\A, f), d=d(\A,f,r)\in \Z.$ These $k$ and $d$  can be easily calculated from the coefficients for $q^2$-terms 
in $a(q,r)$ and $f(p,q)$ by formula (\ref{GaussSum}). The same calculations also prove:

 {\em Either $k(\A, f)=0$ and $\A\s_f[p]$ is a bra-state
or $k(\A,f)\neq 0$ and $\A\s_f[p]$ is pseudo-continuous on $k\Dom+d$ for some $0\le d<k.$}
 
\epk
\bpk A (basic) {\bf Euclidean/Hermitian  Gaussian Hilbert space} $\HH_\Dom$
 is  the $\F$-linear space finitely generated by    Gaussian ket-states.
  and position states. 
 The inner product is defined as the formal Euclidean, respectively, Hermitian inner product between ket-states and ket-states and position states and extends uniquely to their linear combination by bi-linearity law.

\medskip

A {\bf general}  Euclidean/Hermitian 
Gaussian Hilbert space is a tensor power $\HH^{\otimes N}_\Dom$ of  the basic  Gaussian Hilbert space.
\epk
\bpk\label{EH} {\bf Corrsepondence of structures over $\uuu$- \ and $\vvv$-domains}.

Let $\uuu$ and $\vvv$ be as defined in \ref{notation}.

Set
$$\Dom_\vvv:=\vvv\cdot\ZZ[-\frac{\lfr}{2}: \frac{\lfr}{2}-1]\subset \U,
\ \ \NN_\vvv=\lfr$$
and 
$$\Dom_\uuu:=\uuu\cdot\ZZ[-\frac{\lfr\ii}{2}: \frac{\lfr\ii}{2}-1]\subset \U,\ \ \NN_\uuu=\lfr\ii.$$

Note that \be\label{embDom} \Dom_\vvv\subset \Dom_\uuu\mbox{ and }\Dom_\vvv=\ii\cdot \Dom_\uuu.\ee
Indeed, any $\ww\in \Dom_\vvv$ is of the form $\ww=\vvv\cdot l,$ $l\in \ZZ[-\frac{\lfr}{2}:\frac{\lfr}{2}-1] $
and $\vvv=\ii\cdot \uuu.$ That is $$\ww=\uuu \cdot \ii l\in \uuu\cdot \ZZ[-\frac{\lfr\ii}{2}:\frac{\lfr\ii}{2}-1]=\Dom_\uuu.$$ 

Note also that $$\NN_\uuu=\ii \NN_\vvv\mbox{\ and\ }\e(\frac{\ii f(r,p)}{2\NN_\uuu})=\e(\frac{ f(r,p)}{2\NN_\vvv}).$$

The embedding of domains agrees with  the correspondence between the  Gaussian states
$$\{ \}^\ii:\ \s_f[p]\mapsto \s_f^\ii[p]$$
where
$$\s_f[p]=\cc_\s\sum_{r\uuu\in \Dom_\uuu}\e(\frac{f(r,p)}{2\NN_\uuu})\,\uu[r] \mbox{ and }  \s_f^\ii[p]=
 \jj \cc_\s\,\sum_{r\vvv\in \Dom_\vvv}\e(\frac{\ii f(r,p)}{2\NN_\uuu})\,\uu[r]
$$
(recall $\jj=\sqrt{\ii}$).

This can be extended   to the $\F$-linear surjective map
$$\{ \}^\ii: \HH_{\Dom_\uuu}\to  \HH_{\Dom_\vvv}.$$

Consider also the related map
$$\{ \}^\ii: \Oo\GG(\NN_\uuu)\to \Oo\GG(\NN_\vvv); \ \ \ \frac{c}{\sqrt{\NN_\uuu}}\cdot\e(\frac{ f(r,p)}{2\NN_\uuu})\mapsto \frac{c}{\sqrt{\NN_\vvv}}\cdot\e(\frac{ f(r,p)}{2\NN_\vvv})$$

\medskip

Now for formal inner products calculation (\ref{sumInner})  gives:

\be \label{F_inner} \la \s^\ii_{f_1}[p_1]\, |\, \s^\ii_{f_2}[p_2]\ra=\{\la \s_{f_1}[p_1]\, |\, \s_{f_2}[p_2]\ra\}^\ii
\ee
where the inner product is Euclidean on both sides or Hermitian on both sides.


For a linear operator $\A$ of the form (\ref{Ab}) on $\HH_{\Dom_\uuu}$ define the operator  on $\HH_{\Dom_\vvv}$
$$\A^\ii: \  \s^\ii_{f_1}[p]\mapsto \{ A\,\s_{f_1}[p_1]\}^\ii$$

\epk
\bpk \label{Th1} {\bf Theorem.} {\em Under assumptions \ref{EH} the embedding 
$$ \Dom_\vvv\subset \Dom_\uuu$$
gives rise to a canonical surjective homomorphism of Euclidean/Hermitian Gaussian Hilbert spaces  equipped with Gaussian linear maps
$$\{ \}^\ii: \HH_{\Dom_\uuu}\to  \HH_{\Dom_\vvv}.$$ 

The map can be uniquely and definably extended to $M$-dimensional versions of domains    
$$\Dom_\vvv\subset \Dom_\uuu\subset \U^M,\ \ M\in \Oo(\mathcal{F})$$
and to $M$-tensor-power version of Hilbert spaces
$$\{ \}^\ii: \HH_{\Dom_\uuu}^{\otimes M}\to  \HH_{\Dom_\vvv}^{\otimes M}.$$ 
}

{\bf Proof.} The construction and the argument for the first statement is in \ref{EH}. The second statement is just consequence of the algebraic property of finite tensor products. $\Box$
\epk
{
\section{Examples}
\bpk\label{QM} {\bf Example: 1-dimensional QM.}

Set
$$\Dom:=\Dom_\vvv$$
as in \ref{EH}.

Respectively we have position states 
$$\LL_\uu=\{ \uu[r]: r\vvv\in \Dom\}.$$

Define, for $p\vvv\in \Dom$
 the momentum state
$$\vv[p]: r\mapsto \frac{1}{\sqrt{\NN}}\e(-\frac{rp}{\NN})$$

Hence
  $$\vv[p]:=\frac{1}{\sqrt{\NN}}\sum_{r\vvv\in \Dom} \e(-\frac{{rp}}{\NN})\uu[r]$$

 Clearly,
 $$\uu[r]=\frac{1}{\sqrt{\NN}}\sum_{p\vvv\in \Dom} \e(\frac{{rp}}{\NN})\vv[p]$$
and we consider the definable sort
$$ \LL_\vv=\{ \vv[p]: p\vvv\in\Dom \}.$$ 

Thus $\HH_\Dom$ is generated by both the orthogonal systems $\LL_\uu$ and $\LL_\vv$ which are Fourier dual of each other:

One considers the unitary operators $U$ and $V$ that in continuous setting can be written as
 $$U=\e^{i\Qq}\mbox{ and }V= \e^{i\Pp}$$ for the self-adjoint unbounded operators
$\Qq$ (position) and $\Pp$ (momentum).

The operators in our (discrete) setting are defined by their action
$$U: \vv[p]\mapsto \vv[p-{1}]$$
and so acts on $\LL_\uu$ by linearity as
$$U: \uu[r]\mapsto \frac{1}{\sqrt{\NN}}\sum_{p\vvv\in \Dom} \e(\frac{{rp}}{\NN})\vv[p-{1}]=\frac{1}{\sqrt{\NN}}\sum_{p\vvv\in \Dom} \e(\frac{{rp+r}}{\NN})\vv[p]=\e(\frac{{r}}{\NN})\uu[r]$$
and thus the $\uu[r]$ are eigenvectors of the operator.

\medskip

Similarly  the  unitary operator
$$V: \uu[r]\mapsto \uu[r+{1}],\ \ \ \vv[p]\mapsto \e(\frac{{p}}{\NN})\vv[p]$$
has the $\vv[p]\in \LL_\vv$ as its eigenvectors.

It is easy to check that
$$UV=qVU,\mbox{ for } q=\e(\frac{{1}}{\NN})=\exp_\pp(\frac{{1}}{\NN}),\ q^\NN=1.$$

\medskip

{\bf Free particle.} The time evolution operator for the free particle is  $\e^{\frac{i\Pp^2}{2}t},$ where $t\in \ZZ/\NN,$
 the unitary operator
with the action on $\LL_\vv$ defined as 
$$\e^{\frac{i\Pp^2}{2}t} \vv[p]:=\e({\frac{p^2}{2\NN}t})\vv[p].$$

One can calculate 
 $$ \e^{\frac{i\Pp^2}{2}t}\uu[r]=\e^{\frac{i\Pp^2}{2}t} \frac{1}{\sqrt{\NN}}\sum_{p\in \Dom} \e(\frac{rp}{\NN}) \vv[p]=\frac{1}{\sqrt{\NN}}\sum_{p\in \Dom} \e(\frac{p^2t+2rp}{2\NN}) \vv[p]=$$
$$= \frac{1}{\sqrt{\NN}}\sum_{p\in \Dom} \e(\frac{p^2t+2rp}{2\NN})\frac{1}{\sqrt{\NN}}\sum_{s\in \Dom} \e(-\frac{sp}{\NN})\uu[s]=$$ $$=
 \frac{1}{\NN}\sum_{s\in \Dom} \left(\sum_{p\in \Dom}\e(\frac{p^2t+2p(r-s)}{2\NN}) \right)\uu[s]$$
 By  Gauss' quadratic sums formula  
 one gets  
$$\sum_{p\vvv\in \Dom}\e(\frac{p^2t+2p(r-s)}{2\NN}=\left\lbrace\begin{array}{ll}
 \e(-\frac{1}{8})\sqrt{\frac{\NN}{t}} \e(\frac{(r-s)^2}{2t\NN}), \mbox{ if } t|(r-s)\\
0,\mbox{ otherwise}
\end{array}\right. $$
(note that $\e(\frac{1}{8})=\e^\frac{\pi \mathrm{i}}{4}$
 when identifying $\F\subset {^*\C}$).
 
Thus 
\be\label{FreeProp}\e^{\frac{i\Pp^2}{2}t} \uu[r]=\frac{1}{\sqrt{t\NN}}\e(\frac{1}{8}) \sum_{(r-s)\vvv\in t\Dom}
 \e(-\frac{(r-s)^2}{2t\NN})\uu[s] \ee
Note that the range of the state 
$$\mathrm{Range} \left(\e^{\frac{i\Pp^2}{2}t}\right) = t\Dom -r$$
a proper subset of $\Dom,$ a coset of a ``dense'' subgroup.

\medskip

The example of $ \e^{\frac{i\Pp^2}{2}t}$  and similar example $\e^{i\frac{\Pp^2+\omega \Qq^2}{2}t} $ for the {\bf quantum harmonic oscillator} (particles with quadratic potential) can be found in \cite{QM}, sections 11 and 12, and in a more detailed form in \cite{Appro}, section 6.

\epk

\bpk\label{SM} {\bf Example SM.} (Statistical Mechanics) The   setting of SM is even more relevant to the theory above
 since by definition its physics setting is an extremely large but {\bf finite} model $\Dom.$)

In general, models of statistical mechanics are  good analogues of models of {\em quantum field theory}, QFT, rather than quantum mechanics. The similarity becomes apparent once one replaces QFT and QM expressions like $\e^{iS(x)}$ by $\e^{S(x)}$ (Wick rotation). See \cite{McCoy} for a detailed discussion on the topic. 

We single out a more specific SM-setting                                                                                                                                                                                                                                                                                                                                                                                                                                                                                                                                                                                                                                                                                                                                                                                                                                                                                                                                                                                                                                                                                                                                                                                                       
of J.Zinn-Justin, \cite{ZJ}, chapter 4, {\em Classical statistical physics: One dimension}. A more general setting which leads to a real Hilbert space formalism uses {\em probability  density matrices} can also be found in \cite{GusSig}, chapter 9, and indeed in many other sources.

\cite{ZJ} introduces a Hilbert space formalism in SM-context, position and momentum operators (position states but {\bf not} momentum states)   and  
the Gaussian (that is the quadratic case) 
transfer matrix, the propagator between states $q'$ and $q''$
$$T(q',q''):=\e^{S(q',q'')}$$

$T$ a real symmetric operator acting on a real Hilbert space $\HH$ of dimension $\NN.$

                                                                                                                                                                                                                                                                                                                                                                                                                                                                                                                                                                                                                                                                                                                                                                                                                                                                                                                                                                                                                                                                                                                                                                                                       The continuous form of $T(q',q'')$ presented in section 4.6 of \cite{ZJ}                                                                                                                                                                                                                                                                                                                                                                                                                                                                                                                                                                                                                                                                                                                                                                                                                                                                                                                                                                                                                                                                                                                                                                                                                                                                                                                                                                                                                                                                                                                                                                                                                                                                                                                                                                                                                                                                                                                                                                                                                                                                                                                                                                                                                                                                                                                                                                                                                                                                                                                                                                                                                                                                                                                                                                                                                                                                                                                                                                                                                                                                                                                                                                                                                                                                                                                                                                                                                                                                                                                                                           is in full analogy with QM analogue (\ref{exampleHO}) which we discuss below in section \ref{exCL}.

{\tiny                                                                                                                                                                                                                                                                                                                                                                                                                                                                                                                                                                                                                                                                                                                                                                                                                                                                                                                                                                                                                                                                                                                                                                                                        
}

\medskip

We introduce here the domain $\Dom:=\Dom_\mathrm{1SM}:=\Dom_\uuu$ (as in \ref{EH})  of the 1-dimensional statistical mechanics model which agrees with Zinn-Justin's in direct analogy with \ref{QM}
and note that $\NN:=\#\Dom=\lfr\ii=(\mfr\jj)^2$ is not an element of $\F,$ so $\Dom_\mathrm{1SM}$ is not tame.

Now we introduce the ingredients of pseudo-real Gaussian Hilbert space over $\Dom_\mathrm{1SM}.$

Define,  for $r\uuu\in \Dom$ 
$$\uu[r](x)=\delta(x-r).$$

We set, for $p\uuu\in \Dom$ 
$$\vv[p]=\frac{1}{\sqrt{\NN}}\sum_{r\uuu\in \Dom} \e(-\frac{{rp}}{\NN})\uu[r]$$
This makes $\vv[p]$ Fourier-dual to $\uu[r]$  and furnishes two orthonormal bases of $\HH_\Dom$ both of size $\NN=\lfr\ii.$ 
(Note that the above definition of momentum states and its Fourier-duality to position states requires periodicity of $\Dom,$ which in our setting is isomorphic to the group  of period $\lfr\ii,$  much larger scale than $\NN$ of example \ref{QM}.)  

The Gauss quadaratic sums formula appropriate for the given $\Dom$ is presented in (28) of Proposition 4.1 of \cite{FP}, assuming the current system of notations and moving $\mfr$ to LHS
\be\label{G2} \frac{1}{\mfr}\sum_{\eta\in a\Dom} \e(-\frac{a\eta^2 }{2})=\e(\frac{1}{8}) \frac{\jj}{\sqrt{a}}\ee
(bear in mind that $\jj$ specialises by
$\lm_\F$ into $\e(-\frac{1}{8})$  and this brings the RHS to the purely real value $\frac{1}{\sqrt{a}}.$)

\epk
}

\section{Continuous logic setting for Gaussian $\HH_\Dom$.}\label{CL}
\bpk \label{clp} {\bf $\Dom$ in continuous logic CL.
}

We use $\mathsf{lm}_\U: \U\to \bar{\C}$  and $\mathsf{lm}_\F: \F\to \bar{\C}$ 
to move the domain $\Dom$ of states to a locally bounded metric space $\Dom^\lm$ and the domain $\F$ of logical values from $\F$ to $\C.$ 
Our $\Dom$ in this section is either $\Dom_\uuu$ or $\Dom_\vvv$ of \ref{GHDom}.

Continuous logic and model theory is understood as in \cite{HartCMT}.
As in \ref{Dom}
$\Dom$  
is assumed to be given 
with predicates determining non-standard distance between points.

We define metric predicates on the domains assuming that both $\uuu$ and $\vvv$ generate group lattices with the spacing $\frac{1}{\lfr}.$

Thus
$$\mathrm{dist}(0,l\cdot  \uuu):=\frac{l}{\lfr}$$
 $$\mathrm{dist}(0,l\cdot\vvv):=\frac{l}{\lfr}$$
 
 In line with the standard boundedness assumptions  \cite{HartCMT} we distinguish the family of subdomains $\Dom_m\subset \Dom,$ $m\in \N,$
 $$\Dom_m=\{ x\in \Dom:\ \mathrm{dist}(0,x)\le m\}.$$

Set 
 $$\Dom^\lm_m=\mathsf{lm}_\U(\Dom_m)\mbox{ and }\Dom^\lm=\bigcup_{m\in \N} \Dom^\lm_m.$$
Thus
$\Dom^\lm $ is covered by a family $\Dom^\lm_m\subseteq \Dom^\lm$  of bounded complete metric subspaces 
$\Dom^\lm_m.$ Note that in general $\Dom\neq \bigcup_{m\in \N} \Dom_m.$

 It follows $$\Dom^\lm_\uuu\cong \Dom^\lm_\vvv\cong \R$$ as metric spaces and
 $\Dom_m^\lm$ in both cases corresponds to $\R\cap [-m: m].$ However, note that since by our choices of parameters 
 $$\Dom_\vvv=\ii \Dom_\uuu\mbox{ and } \lm_\U (\ii u)= \mathrm{i}\, \lm_\U(u),$$
 more accurately, we have isomorphisms 
 \be\label{DomcGauss} \Dom^\lm_\uuu\cong \R\mbox{ and }\Dom^\lm_\vvv\cong \mathrm{i}\R\ \ \ (\mathrm{i}=\sqrt{-1})\ee

 We use the fact that {\em $\Dom$ has a structure of a $K$-module}, for a   subring $K\subseteq \ZZ$ of our choice.
 For most purposes we can consider $K=\Z$ or slightly bigger subring of $\Oo(\mathcal{F}$). For $k\in K,$ $k\Dom$ is a submodule and, for $r\in \Dom,$ $k\Dom+ r$ is a coset of the submodule,
\be\label{sizekDom} \#(k\Dom+r)=\frac{\NN}{k}\ee
and at the same time 
\be \label{lmkDom} \mathsf{lm}_\U\, (k\Dom+r)=\Dom^\lm\ee
(we say that $k\Dom$ is dense in $\Dom$).

Let $\varphi: \Dom\to \F$ be a pseudo-continuous definable state on $\Dom.$ Set for $x\in \Dom_m,$   
\be\label{ket}x^\lm:=\mathsf{lm}_\U(x)\mbox{ and }\varphi^\lm(x^\lm):= \lm_\F(\sqrt{\NN}\varphi(x)).\ee

This is well-defined because of pseudo-continuity and defines a map
$$\varphi^\lm: \Dom^\lm \to \C.$$

In the more general case we consider $\varphi$ which is pseudo-continuous on $k\Dom+r$ as above and apply the definition (\ref{ket}) to
$\varphi_{|k},$ the restriction to $k\Dom+r.$ This is well-defined.

\epk

\bpk\label{7.3} {\bf Lemma.} {\em Let $\varphi$ be a ket state on $\Dom$ of norm 1  and 
$\varphi_{|k}$ its restriction on $k\Dom$ normalised to norm 1. Then
\be\label{kdom}\varphi_{|k}^\lm(x^\lm)=\varphi^\lm(x^\lm)\mbox{ for }x\in k\Dom+r.\ee

More generally, the equality holds for any ket state $\varphi$ on $\Dom$ 
when setting $\varphi_{|k}: (k\Dom+r) \to \F$ 
$$\varphi_{|k}(x):=\sqrt{k}\cdot\varphi_{|k}(x).$$
} 

{\bf Proof.} We may assume $r=0.$

By definition $$\varphi=\sum_{s\in \Dom} \varphi(s)\uu[s],$$
$$\varphi_{|k}=\sqrt{k}\sum_{s\in k\Dom} \varphi(s)\uu[s]$$ and so
$$\varphi_{|k}(s)=\sqrt{k}\varphi(s).$$

Thus by (\ref{ket}), for $s\in k\Dom,$
$$\varphi^\lm_{|k}(s^\lm)=\mathsf{lm}\left(\sqrt{\frac{\NN}{k}}\sqrt{k}\varphi(s)\right)=\varphi^\lm(s^\lm).$$
$\Box$

\epk

The CL-state of the form $\varphi^\lm_{|k}$ we call CL-ket states. Lemma \ref{7.3} allows us always to make reductions, if needed, to appropriate dense subdomains. 

Define {\bf  inner product for CL-ket-states} on $\Dom^\lm,$ for $\varphi,\psi,$ $\varphi\neq \psi$
\be \label{defInner}  \la\varphi^\lm_\F|\psi^\lm\ra=\lim_{m\to \infty} \int_{\Dom^\lm_m} \varphi^\lm(x){\bar{\psi}}^\lm(x)\, dx
\ee
If $\varphi= \psi$ then set 
\be\label{37'}\la\psi^\lm|\psi^\lm\ra:=\lm_\F |\psi|^2.\ee

\bpk{\bf Proposition.}   $\la\varphi^\lm|\psi^\lm\ra$ has a well-defined value in $\C.$ More precisely 
\be \label{innerket}\la\varphi^\lm|\psi^\lm\ra=\left\lbrace\begin{array}{ll}
\lm_\F (\mfr \la\varphi|\psi\ra_\mathrm{E}), \mbox{ if } \la\varphi|\psi\ra_\mathrm{E}\neq 0\ \&\  \Dom=\Dom_\uuu\\
\lm_\F (\mfr \la\varphi|\psi\ra_\mathrm{H}), \mbox{ if } \la\varphi|\psi\ra_\mathrm{H}\neq 0\ \& \ \Dom=\Dom_\vvv
 \end{array}\right. 
\ee

{\bf Proof.} We assume without loss of generality that $|\varphi|=1=|\psi|.$
Then by assumptions we may write
$$\varphi(r)=\frac{1}{\sqrt{\NN}}\e(-\frac{A_\phi r^2+2B_\phi r}{2\NN})\mbox{ and }\psi(r)=\frac{1}{\sqrt{\NN}}\e(-\frac{A_\psi r^2+2B_\psi\cdot r}{2\NN})$$
$A_\phi\ge  0,\ A_\psi\ge 0.$

Now assume that $\Dom=\Dom_\uuu$ and $\la\varphi|\psi\ra_\mathrm{E}\neq 0$
($\NN=\NN_\uuu=\lfr\ii$). The second assumption implies by formula (\ref{GaussSum}) that $A:=A_\phi+A_\psi>0,$ $B=B_\phi+B_\psi$ is divisible by $A$ and  
$$\la\varphi|\psi\ra_\mathrm{E}=\frac{1}{\NN}\sum_{rA\uuu\in \Dom}
\e(- \frac{A r^2+2Br}{2\NN})=\frac{1}{\sqrt{\NN}}\e(-\frac{1}{8})\sqrt{\frac{1}{A}}\e(-\frac{\frac{B^2}{2A}}{\NN})$$
and so, since $\sqrt{\NN}=\mfr\jj,$
$$\mfr\cdot \la\varphi|\psi\ra_\mathrm{E}=\jj\inv \e(-\frac{1}{8}) \sqrt{\frac{1}{A}}\e(-\frac{\frac{B^2}{A}}{\NN})=\jj\inv \e(-\frac{1}{8})\sqrt{\frac{1}{A}}\exp_\pp(-\frac{\frac{B^2}{A}(\pp-1)}{ \lfr\ii})=$$
$$=\jj\inv \e(-\frac{1}{8})\sqrt{\frac{1}{A}}\exp_\pp(-\frac{B^2}{A} \uuu)$$
Hence, since by definition (see \ref{defIF} - \ref{3.8})
$$\lm_\F:\ \ \jj\mapsto   \e^{i\frac{\pi}{4}},\ \ \e(-\frac{1}{8})\mapsto \e^{i\frac{\pi}{4}},\ \ \exp_\pp(-\frac{B^2}{A} \uuu)\mapsto 1,$$
we obtain
$$\lm_\F\  \mfr\cdot \la\varphi|\psi\ra_\mathrm{E}= \sqrt{\frac{1}{A}}$$


On the other hand,
by definition, for $x=r^\lm,$
 $$\varphi^\lm(x)=\e^{-\pi (A_\phi x^2+2B_\phi x)}
 \mbox{ and }\psi^\lm(x)=\e^{-\pi(A_\psi x^2+2B_\psi x)}$$
 and so we have the convergence for the integral
 $$\int_{x\in \R} \varphi^\lm\cdot \bar{\psi}^\lm \, dx= \int_{x\in \R} \e^{-\pi(Ax^2+2Bx)} dx=\sqrt{\frac{1}{A}}$$
 the classical Gauss integral. This proves the Euclidean case.
 
 \medskip
 
 Now assume that $\Dom=\Dom_\vvv$ and $\la\varphi|\psi\ra_\mathrm{H}\neq 0$
($\NN=\NN_\vvv=\lfr$). The second assumption implies by formula (\ref{GaussSum}) that $A:=A_\phi-A_\psi\neq 0,$ $B=B_\phi-B_\psi$ is divisible by $2A$ and  (assuming for simplicity $A>0$)
$$\la\varphi|\psi\ra_\mathrm{H}=\frac{1}{\NN}\sum_{rA\uuu\in \Dom}
\e(- \frac{A r^2+Br}{2\NN})=\frac{1}{\sqrt{\NN}}\e(-\frac{1}{8})\sqrt{\frac{1}{A}}\e(-\frac{B^2}{2A\NN})$$
and so, since $\sqrt{\NN}=\mfr,$
$$\mfr\cdot \la\varphi|\psi\ra_\mathrm{H}=\e(-\frac{1}{8}) \sqrt{\frac{1}{A}}\e(-\frac{B^2}{2A\NN})= \e(-\frac{1}{8})\sqrt{\frac{1}{A}}\exp_\pp(-\frac{B^2}{2A}\frac{\pp-1}{ \lfr})=$$
$$= \e(-\frac{1}{8})\sqrt{\frac{1}{A}}\exp_\pp(-\frac{B^2}{2A} \vvv)$$
Since 
$$\lm_\F:\ \   \exp_\pp(-\frac{B^2}{2A} \vvv)\mapsto 1$$
we obtain
$$\lm_\F\  \mfr\cdot \la\varphi|\psi\ra_\mathrm{H}= \e^\frac{\pi i}{4}\sqrt{\frac{1}{A}}$$


On the other hand,
by definition, for $x=r^\lm,$
 $$\varphi^\lm(x)=\e^{-\pi i(A_\phi x^2+B_\phi x)}
 \mbox{ and }\psi^\lm(x)=\e^{-\pi i(A_\psi x^2+B_\psi x)}$$
 and so we see, for each $m$, the Fresnel integral
 $$\int_{x\in \Dom^\lm_m} \varphi^\lm\cdot \bar{\psi}^\lm dx= \int_{-m\le x\le m} \e^{-\pi i(Ax^2+Bx)}dx$$
 which has well-defined limit as $m\to \infty,$
 $$\lim_{m\to \infty} \int_{-m\le x\le m} \e^{-\pi i(Ax^2+Bx)}dx=\e^\frac{\pi i}{4}\sqrt{\frac{1}{A}}$$
 
  This proves the Hermitian case and finishes the proof. $\Box$
 
\epk



\bpk \label{7.5}
There will be a different treatment for position states $\uu[r]$
$$\la \varphi^\lm\,|\,\uu^\lm[r^\lm]\ra: =\varphi^\lm(r^\lm)$$
See also (\ref{ket}).

\epk

 \bpk Define for ket-sorts 
 $$\LL^\lm_f:=\{ \s_f^\lm[p]: \ p\in \Dom^\lm\}$$
  
For $x\in \Dom^\lm$
$$ \uu^\lm[x]: p\mapsto \la \vv^\lm[p]\,|\, \uu^\lm[x]\ra=\e^{ipx};\ \ \LL^\lm_\vv\to \C$$
is continuous in $p$ and thus we may identify $\uu^\lm[x]$ with an element of the space of continuous maps
$$\LL^\lm_\vv\to \C,$$
equivalently of continuos linear functionals on $\HH_{\Dom^\lm}.$

\medskip

Using (\ref{defInner}) and (\ref{37'}) we extend  the definition of inner  product to finite linear combinations of elements of sorts $\LL_\phi$ for all $\phi.$
\medskip

{\bf Remark.} Recall that  bi-linear operators on the space of compactly supported functions can be represented by integrals and compactly supported functions are dense in the space of continuous functions on a locally compact set, the definition (\ref{innerket}) is in physics Dirac calculus is represented as
\be \label{innerprodInt} \la f|g\ra=  \int_{\Dom^\lm} f(x)\bar{g}(x) dx 
\ee 
where $f(x)=\phi^\lm(x),\ g(x)=\psi^\lm(x)$ and the integration is over the measure given by Dirac delta-function $\delta_\mathrm{Dir}.$

In particular,
\be\label{dirac?} \frac{1}{2\pi}\int_{\Dom^\lm} \e^{i(p_1-p_2)x} dx=\delta_\mathrm{Dir}(p_1-p_2)\ee
A rigorous theory of Dirac integration (\ref{innerprodInt}) can be found e.g. in \cite{ZeidlerI}, section 10.2. 
\epk

\bpk
The  {\bf ket-CL-Hilbert space} $\HH_{\Dom^\lm}$ is presented as a  multisorted structure with sorts  $\LL_\psi$ 
 for all ket $\psi$
and binary CL-predicates (Dirac inner product)
$$\la \varphi|\psi\ra:\  \LL_\varphi\times \LL_\psi\to \C.$$
defined in (\ref{innerket}).

\medskip
A {\bf bra} CL-state $\psi^\lm$ is the linear functional 
$$\psi^\lm: \LL^\lm_\varphi \to \C: \ \mbox{ for }y\in \Dom^\lm,\  \psi^\lm(y):=\lm \sqrt{\NN}\la \psi| \varphi_y\ra $$ 

This is well-defined by \ref{inn} and coincides with the product on kets when the bra-state $\psi^\lm$ originates from a ket $\psi.$ 

\epk
\bpk\label{7.6}
Define the (full)  {\bf Hermitian CL-Hilbert space} $\HH_{\Dom^\lm}$ as 
  multisorted structure with sorts  $\LL^\lm_\psi$ where $\psi$ can be both bra- or ket-, with Dirac inner product as  binary CL-relation between sorts.

In this treatment of Hilbert-space formalism it is natural to consider  
 $\HH^\mathrm{eq}_{\Dom^\lm}$ which contains finite linear combinations of elements of $\HH_{\Dom^\lm}.$

\epk
\bpk\label{7.8} Linear unitary operators on $\HH_{\Dom^\lm}$ have form
$$\A^\lm: \LL^\lm_\psi\to \LL^\lm_\varphi$$
where $$A:\LL_\psi\to \LL_\varphi$$ is an operator on $\HH_\Dom$ as described in \ref{operators} and \ref{product} as

 $$A: \psi_r\mapsto \sum_{r\in \Dom}\left(\sum_{s\in \Dom}\psi(r,s)\cdot \alpha(s,r)\right)\uu[r]$$
 
 and 
 $$\A^\lm: \psi^\lm_r\mapsto \left(\sum_{r\in \Dom}\left(\sum_{s\in \Dom} \psi(r,s)\cdot\alpha(s,r))\right)\uu[r]\right)^\lm.$$
  
Note that  
 $$r\mapsto \sum_{s\in \Dom} \psi(r,s)\cdot\alpha(s,r)=\la \psi_r|\bar{\alpha}_r\ra$$
 and thus the integral expression
  (\ref{innerprodInt}) for inner product 
brings us   
in CL-setting to
\be\label{25} \A^\lm: \psi^\lm_r\mapsto \int_{\Dom^\lm} \psi^\lm(r,s)\cdot\alpha^\lm(s,r)\, ds\ee
the Dirac integral as a {\bf CL-quantifier} bounded by condition $\alpha$.  
\epk

{
\bpk{\bf Example. One-dimensional quantum mechanics}\label{exCL}

We follow \ref{QM}. Note that we assume $\hbar=1$ and suppress some normalising coefficients. 
 
 A momentum CL-state (ket-state) is
$$\vv[p]: x\mapsto \e^{i p x},\ \mbox{for }x,p\in \R$$


Respectively
$$\LL_\vv=\{ \vv[p]:\ \ p\in \R\}$$

 A position state is a bra-state
a map (that is a CL-unary predicate):  $$\uu[r]: \ \LL_\vv\to \Ss\subset \C;\ \ \vv[p]\mapsto \la \vv[p]|\uu[r]\ra:= \e^{ipr}$$
and the sort 
$\LL_\uu$ with the pairing (binary CL-predicate)
$$ \la \vv[p]|\uu[r]\ra: \ \LL_\vv\times \LL_\uu\to \Ss\subset \C.$$


Both position and momentum  are binary CL-predicates on $\R$ as are  general CL-states.


\medskip


\medskip


\epk

   

\bpk In \cite{QM} and \cite{Appro} a number of calculations in Gaussian setting were carried out.
In  \cite{Appro} we calculated the propagator for the quantum harmonic oscillator
 with frequency $\omega,$ the CL-value of reaching the position $x$ from position $x_0$  in time $t$: 
\be\label{exampleHO}\e^{-i\frac{\pi}{4}} \sqrt{\frac{\omega}{2\pi \hbar |\sin \omega t|}}\exp \ i\omega \frac{(x^2+x_0^2)\cos \omega t-2xx_0}{2\hbar\sin \omega t}\ee 
 
 This was demonstrated therein both by path integral calculation (section 9) and by a more direct calculation (section 7,  7.13). 
 
\epk
\bpk {\bf Conclusions.} {\em Equalities of the form 
 $$\A^\lm ( \psi^\lm) =\varphi^\lm,$$
 for Gaussian $\psi, \varphi$ and $\A,$
 obtained as the result of the calculus defined by (\ref{defInner}) - (\ref{25}) are CL-sentences which form a CL-theory with the interpretation by pseudo-finite structures based on $\Dom_\vvv$ or $\Dom_\uuu.$ 
 
 This theory has quantifier-elimination to Gaussian predicates (states) since by Gaussian summation formula the application of an $\A$ to a $\psi$ always results in a Gaussian state $\varphi.$
 
Wick rotation $$\{ \}^\ii: \HH_{\Dom_\uuu}\to  \HH_{\Dom_\vvv}$$
described by \ref{Th1} establishes an equaivalence between the theories 
 based on $\Dom_\vvv$ or $\Dom_\uuu.$ 
   }
\epk 

\end{document}